\documentclass[11pt,leqno]{amsart}

\usepackage[cp1252]{inputenc}
\usepackage[english]{babel}
\usepackage[a4paper]{geometry}
\usepackage{amsmath,amsfonts,amssymb,amsthm,cases,upgreek}
\usepackage{empheq}
\usepackage{stmaryrd}
\usepackage{fullpage}

\usepackage{diagbox}

\usepackage{dsfont}
\usepackage{graphicx}
\usepackage{comment}
\usepackage{bm} 
\usepackage{dsfont}

\newcommand\iD{\textrm{i}}

\newcommand{\sgn}{\text{\rm sgn}}

\newcommand\br{\begin{remark}}
\newcommand\er{\end{remark}}
\newcommand\bp{\begin{pmatrix}}
\newcommand\ep{\end{pmatrix}}
\newcommand{\be}{\begin{equation}}
\newcommand{\ee}{\end{equation}}
\newcommand{\bes}{\begin{equation*}}
\newcommand{\ees}{\end{equation*}}
\newcommand\ba{\begin{equation}\begin{aligned}}
\newcommand\ea{\end{aligned}\end{equation}}
\newcommand\bas{\begin{equation*}\begin{aligned}}
\newcommand\eas{\end{aligned}\end{equation*}}

\newcommand{\beg}{\begin{example}}
\newcommand{\eeg}{\end{exaplem}}
\newcommand{\bpr}{\begin{proposition}}
\newcommand{\epr}{\end{proposition}}
\newcommand{\bt}{\begin{theorem}}
\newcommand{\et}{\end{theorem}}
\newcommand{\bc}{\begin{corollary}}
\newcommand{\ec}{\end{corollary}}
\newcommand{\bl}{\begin{lemma}}
\newcommand{\el}{\end{lemma}}
\newcommand{\bd}{\begin{definition}}
\newcommand{\ed}{\end{definition}}
\newcommand{\brs}{\begin{remarks}}
\newcommand{\ers}{\end{remarks}}

\newtheorem{theorem}{Theorem}[section]
\newtheorem{proposition}[theorem]{Proposition}
\newtheorem{corollary}[theorem]{Corollary}
\newtheorem{lemma}[theorem]{Lemma}

\theoremstyle{remark}
\newtheorem{remark}[theorem]{Remark}
\theoremstyle{definition}
\newtheorem{definition}[theorem]{Definition}

\newtheorem{example}[theorem]{Example}

\newcommand\R{\mathbb R}

\newcommand{\N}{\mathbb N}
\newcommand{\Z}{\mathbb Z}
\newcommand\bU{{\mathbf U}}
\newcommand\bV{{\mathbf V}}
\newcommand\bW{{\mathbf W}}

\usepackage{xcolor}
\newcommand{\colr}{\color{black}}
\newcommand{\colb}{\color{black}}
\newcommand{\colp}{\color{black}}

\def \epsilon {\varepsilon}

\SetSymbolFont{stmry}{bold}{U}{stmry}{m}{n}

\usepackage{url}
\usepackage{textcase}
\usepackage{hyperref}

\title{\textbf{Rigid lid limit in shallow water over a flat bottom}}

\author{Benjamin Melinand}
\address{CEREMADE, CNRS, Universit\'e Paris-Dauphine, Universit\'e PSL, 75016 Paris, France}
\email{{\tt melinand@ceremade.dauphine.fr}}
\thanks{}

\begin{document}

\begin{abstract}
\colr \noindent We perform the so-called rigid lid limit on different shallow water models such as the abcd Bousssinesq systems or the Green-Naghdi equations. To do so we consider an appropriate nondimensionalization of these models where two small parameters are involved: the shallowness parameter $\mu$ and a parameter $\epsilon$ which can be interpreted as a Froude number. When parameter $\epsilon$ tends to zero, the surface deformation formally goes to the rest state, hence the name rigid lid limit. We carefully study this limit for different topologies. We also provide rates of convergence with respect to $\epsilon$ and a careful attention is given to the dependence on the shallowness parameter $\mu$. \colb
\end{abstract}

\date{\today}
\maketitle

\section{Introduction}

We consider shallow water asymptotic models of the water wave equations that can be written under the general form
\[
\epsilon \partial_{t} \bU + \mathcal{L}_{\mu} \bU + \epsilon Q_{\mu}(\bU) = 0.
\]
Here $t \in \R^{+}$, \colr $x \in \R$ or $\R^{2}$ \colb, $\bU=(\zeta,\bV) \in \R^{2}$ or $\R^{3}$, $1+\epsilon \zeta$ is the nondimensionalized water depth, $\bV$ is the nondimensionalized horizontal velocity of the fluid, $\mathcal{L}_{\mu}$ is a linear operator that tends to the wave operator as $\mu \to 0$ and $Q_{\mu}$ a nonlinear operator. The nondimensionalized parameter $\mu$ measures the shallowness of the flow. The nondimensionalized parameter $\epsilon$ compares the amplitude of the water waves to the water depth and can also be seen as a comparison between the typical horizontal velocity of the fluid and the typical velocity of the water waves. Note that the nondimensionalized parameter $\epsilon$ appears in front of $\partial_{t}$ since the characteristic time used to nondimensionalize the time variable is the typical time scale of the fluid (and not the one of water waves). We assume in the following that $\epsilon$ and $\mu$ both belong to $(0,1]$. 

The goal of this paper is to perform the limit $\epsilon \to 0$ and understand how the parameter $\mu$ interferes with this convergence. \colr In our framework the velocity component $\bV$ tends to a solution of the incompressible Euler equations whereas \colb the free surface $\zeta$ tends to $0$ thus the terminology rigid lid limit. Actually the rigid lid approximation is a common assumption in the oceanographical literature and can be interpreted as a low Froude number assumption (ocean currents travel slower than water waves). We refer for instance to \cite{Bowen_Holman89,dodd_thornton_90,allen_newberger_holman_96} where existence and stability of nearshore shear waves that are not gravity waves were discussed and such assumption is used. In the same direction Camassa, Holm and Levermore  \cite{CHL96,CHL97} derive from the free surface Euler equations two asymptotic models called nowadays the lake equations and the great lake equations where again a rigid lid approximation is considered. \colp In our setting the lake equations are obtained by neglecting terms of order $\mathcal{O}(\mu)$ and then performing the limit $\epsilon \to 0$ whereas the great lake equations are derived by neglecting terms of order $\mathcal{O}(\mu^2)$ and then performing the limit $\epsilon \to 0$\colb. A full justification of the lake equations was obtained in \cite{bresch_metivier} (see also \cite{oliver97}). We emphasize that all the previous cited works consider a non flat bathymetry. That is not the case in this paper where the seabed is assume to be flat. Finally we also refer to \cite{Benoit_rigid_limit} where a first study of the rigid lid limit on the full irrotational water wave equations is performed. The proof is however based on weighted dispersive estimates that are not well suited for the local wellposedness on large time and the \colr rates of convergence \colb obtained are not optimal. \colr Using the strategy provided in this paper one can improve \cite{Benoit_rigid_limit} by obtaining rates of convergence as those established for instance in Section \ref{s:GN} for the irrotational Green-Naghdi equations.\colb 

We carefully study if the convergence is strong, meaning in $L^{\infty}(0,T;L^{2}(\R^{n}))$ for some $T>0$ independent of $\epsilon$. If not, we provide the default of strong convergence, meaning a corrector $\hat{\bU}_{\epsilon}$ such that $\bU-\hat{\bU}_{\epsilon}$ strongly converges. We also get convergence in $L^{q}(0,T;L^{\infty}(\R^{n}))$ for some $q \geq 2$ (weak convergence) and in $L^{2}(0,T;L^{2}_{loc}(\R^{n}))$ (convergence of the local energy)\footnote{\colr Note that this kind of convergence was also obtained for the rigid lid limit of the shallow water equations over a nonflat bottom in \cite{bresch_metivier} (without any rate of convergence), see Theorem 2.5 and Subsection 4.2.2 inside.\colb}. All our convergence results are given with \colr rates of convergence \colb in terms of $\epsilon$ that crucially depends on the shallowness parameter $\mu$ and the function space we consider. \colr Our strategy is the following. We assume bounds on appropriate $H^s$-norms of $\bU$ that are uniform with respect to $\epsilon$ on a existence time\footnote{\colr Such bounds and existence time can be obtained from energy estimates when one prove the local wellposedness, see for instance Subsections \ref{localexist1d} or \ref{localexist2d} below.\colb} that is independent of $\epsilon$. Using Duhamel's formulation one can then use Strichartz estimates to get $L^{q}_{t}L^{r}_{x}$ bounds and Morawetz-type estimates to obtain $L^{2}_{t} L^{2}_{loc}$ controls.\colb

Our problem shares some similarities with the incompressible limit (low Mach number limit) of the compressible Euler equations defined in $\R^{n}$: 
\bes
\left\{
\begin{aligned}
&\epsilon \partial_{t} c + \epsilon \bV \cdot \nabla  c + \frac{\gamma-1}{2} (1+\epsilon c) \nabla \cdot \bV = 0,\\
&\epsilon \partial_{t} \bV + \epsilon (\bV \cdot \nabla) \bV + \frac{\gamma-1}{2} (1+\epsilon c) \nabla c = 0.
\end{aligned}
\right.
\ees
In that case $1+\epsilon c$ is the rescaled speed of the sound, $\bV$ is the velocity of the fluid, $\gamma>1$ the adiabatic exponent and $\epsilon$ the Mach number. The limit $\epsilon \to 0$ was studied by several authors and we refer for instance to \cite{Ukai86,MetivierSchochet01,DutrifoyHmidi04,Gallagher05,Alazard08,HmidiSulaiman14}. It is now well understood that the acoustic component of $(c,\bV)$ is of dispersive type and weakly converges to $0$ since its propagation speed is of size $\frac{1}{\epsilon}$ and, when the initial datum only contains incompressible terms, the convergence is strong. Such phenomenon also appears for the rigid lid limit.

\colp Finally we mention the following works \cite{LPS12,SautXu20,SautXu21,Tesfahun24} where dispersive estimates similar to ours are used to study the long time existence problem on Boussinesq and Boussinesq-type systems.\colb

We organize the paper as follows. Section \ref{s:classicalBoussi1d} and Section \ref{s:classicalBoussi2d} are devoted to the study of the rigid lid limit on the classical Boussinesq system (also called Amick-Schonbek system in the literature) respectively in 1d and 2d. We explain in details our strategy. In Section \ref{s:Boussiabcd} we consider other Boussinesq systems. Finally in Section \ref{s:GN} we perform the rigid lid limit on the Green-Naghdi equations. Appendix \ref{s:LP} recall some basic facts about Littlewood-Paley decomposition and Appendix \ref{s:fouriermulti} gathers different useful Fourier multiplier estimates. Finally we provide in Appendix \ref{s:dispestim} general dispersive estimates on linear dispersive equations with radial nonhomogeneous phases. 

\subsection*{Notations}
\begin{itemize}
\item If $f$ is a Schwartz class function defined on $\R^{n}$, we define $\mathcal{F}f$ or $\hat{f}$ as the Fourier transform of $f$ by
\[
\hat{f}(\xi) = \frac{1}{(2\pi)^{\frac{n}{2}}} \int_{\R^{n}} e^{-\iD x \cdot \xi} f(x) dx.
\]

\item \colp If $m$ is a smooth function that is at most polynomial at infinity, we define the Fourier multiplier $m(D)$ as, for any Schwartz class function $f$, 
\[
m(D) f = \mathcal{F}^{-1} (  m(\xi) \hat{f}(\xi)).
\]\colb
Note that the Laplace operator verifies $\Delta=-|D|^{2}$.

\item If $f$ is function defined on $\R^{n}$ we denote by $\| f \|_{L^{p}}$ the $L^{p}(\R^{n})$ norm of $f$.

\item \colr If $p \in [1,\infty]$, we denote $p'=\frac{p}{p-1}$.\colb

\item If $T>0$ and $G : (t,x) \in [0,T] \times \R^{n} \to G(t,x) \in \R$, the norm $\| G \|_{L^{q}_{t} L^{r}_{x}}$ corresponds to the norm of the space $L^{q}((0,T);L^{r}(\R^{n}))$. 

\item If $S : E \to F$ is a linear bounded operator with $E,F$ two Banach spaces, we denote by $S^{\ast}$ its adjoint.
\end{itemize}

\subsection*{Acknowledgments} \colr We deeply thank the anonymous referee for careful reading and for very valuable comments on the manuscript. \colb This work has been partially funded by the ANR project CRISIS (ANR-20-CE40-0020-01).

\section{The classical 1D Boussinesq system}\label{s:classicalBoussi1d}

In the one dimensional case the classical Boussinesq system reads as

\be\label{boussinesq_eq1d}
\left\{
\begin{aligned}
&\epsilon \partial_{t} \zeta + \partial_{x} \left( \left[1 + \epsilon \zeta \right] V \right) = 0,\\
&\epsilon \left(1- \frac{\mu}{3} \partial_{x}^{2} \right) \partial_{t} V + \partial_{x} \zeta + \epsilon V \partial_{x} V = 0.
\end{aligned}
\right.
\ee
Denoting by $\bU = (\zeta, V)^{T} \in \R^{2}$, we get the following system
\be\label{boussinesq_eq1dcompactform}
\epsilon \partial_{t} \bU + A(\partial_{x}) \bU = \epsilon F(\bU)
\ee
where
\bes
A(\partial_{x}) = \begin{pmatrix} 0 & \partial_{x} \\ (1- \frac{\mu}{3} \partial_{x}^{2})^{-1} \partial_{x} & 0 \end{pmatrix} \text{   and   } F(\bU) = - \begin{pmatrix} \partial_{x} (\zeta V) \\ \frac12 (1- \frac{\mu}{3} \partial_{x}^{2})^{-1} \partial_{x} (V^2) \end{pmatrix}. 
\ees
Note that
\be\label{exp(tA)1d}
\exp(tA(\partial_{x})) = \begin{pmatrix} \cos \left( t \frac{D}{\sqrt{1 + \frac{\mu}{3} D^2}} \right) & \iD \sqrt{1 + \frac{\mu}{3} D^2} \sin \left( t \frac{D}{\sqrt{1 + \frac{\mu}{3} D^2}} \right) \\ \frac{\iD}{\sqrt{1 + \frac{\mu}{3} D^2}} \sin \left( t \frac{D}{\sqrt{1 + \frac{\mu}{3} D^2}} \right) & \cos \left( t \frac{D}{\sqrt{1 + \frac{\mu}{3} D^2}} \right) \end{pmatrix}.
\ee

We state the main result of this section.

\bt\label{rigidlid1d}
Let $M>0$, $T>0$, $\epsilon \in (0,1]$ and $\mu \in (0,1]$. Let $(\zeta,V) \in \mathcal{C}([0,T];(H^{3} \times H^{3})(\R))$ a solution of \eqref{boussinesq_eq1d} with initial datum $(\zeta_{0},V_{0})$ such that
\[
\| (\zeta,V) \|_{L^{\infty}(0,T;H^{3}(\R) \times H^{3}(\R))} \leq M.
\]
There exists a constant $C>0$ \colr independent of $M$, $T$, $\epsilon$ and $\mu$ \colb such that for any $q,r \geq 2$ with $\frac{1}{q} + \frac{1}{2r} = \frac{1}{4}$ and any $\tilde q, \tilde r \geq 2$ with $\frac{1}{\tilde q} + \frac{1}{3 \tilde r} = \frac{1}{6}$
\[
\begin{aligned}
&\left\| \bp \zeta \\ V \ep - e^{-\frac{t}{\epsilon} A(\partial_{x})} \bp \zeta_{0} \\ V_{0} \ep \right\|_{L_{t}^{q}(0,T;L_{x}^{r}(\R))} \leq \left( \frac{\epsilon}{\mu} \right)^{\frac14 + \frac{1}{q}} \colr M^2 \colb  T^{\frac34} C,\\
&\colp \left\| \bp \zeta \\ V \ep - e^{-\frac{t}{\epsilon} A(\partial_{x})} \bp \zeta_{0} \\ V_{0} \ep \right\|_{L_{t}^{\tilde q}(0,T;L_{x}^{\tilde r}(\R))} \leq \left( \frac{\epsilon}{\mu} \right)^{\frac16 + \frac{1}{\tilde q}}  M^2  T^{\frac56} C, \colb\\ 
&\left\| \bp \zeta \\ V \ep \right\|_{L_{t}^{\tilde q}(0,T;L_{x}^{\tilde r}(\R))} \leq \left( \frac{\epsilon}{\mu} \right)^{\frac{1}{\tilde q}} \left( \colr M \colb + \colr M^2 \colb T^{\frac56} \left( \frac{\epsilon}{\mu} \right)^{\frac{1}{6}} \right) C,\\
&\sup_{x_{0} \in \R} \left\| e^{-(x-x_{0})^2} \bp \zeta \\ V \ep \right\|_{L_{t}^{2}(0,T;L^{2}_{x}(\R))} \leq \epsilon^{\frac12}(\colr M \colb + \colr M^2 \colb T) C.
\end{aligned}
\]
\et 

\br 
Several remarks are in order. We first note that $(\zeta,V)$ weakly converges to $0$ as $\epsilon \to 0$. However $(\zeta,V)$ does not strongly converge to $0$, meaning in $L^{\infty}(0,T;L^{2}(\R))$, except if $(\zeta_{0} , V_{0})=(0,0)$. The first two estimates provides the default of strong convergence. The first three estimates are not uniform with respect to $\mu \to 0$. The fourth estimate provides the convergence to $0$ in $L^{2}_{t}(0,T;H^{s}_{loc}(\R))$ for any $s \in [0,3)$ and is uniform with respect to $\mu \to 0$. That corresponds to the decay to $0$ of the local energy. 
\er

\br\label{remarkKdV} Note that if $\epsilon \sim \mu$ as in \cite{bona_chen_saut_derivation,bona_colin_lannes} or when \colr$\mu = \mathcal{O}(\epsilon)$\colb, the first two estimates does not provide a convergence result as $\epsilon \to 0$. It is known in that case that nonlinear terms must be taken into account and that asymptotic models like a system of decoupling KdV equations or a system of decoupling BBM equations become relevant. 
\begin{comment}
To be more precise, if we are under the assumptions of Theorem \ref{rigidlid1d} and if $g_{\pm}$ satisfy
\[
\left\{
\begin{aligned}
&\epsilon \partial_{t} g_{+} + \partial_{x} g_{+} + \frac{\mu}{6} \partial_{x}^{3} g_{+} + \frac32 \epsilon g_{+} \partial_{x} g_{+} = 0,\\
&\epsilon \partial_{t} g_{-} - \partial_{x} g_{-} - \frac{\mu}{6} \partial_{x}^{3} g_{-} - \frac32 \epsilon g_{-} \partial_{x} g_{-} = 0,\\
&g_{+}(0,x)= \frac12 (\zeta_{0}(x) + V_{0}(x)) \text{   ,   } g_{-}(0,x)=\frac12(\zeta_{0}(x) - V_{0}(x)),
\end{aligned}
\right.
\]
there exists a constant $C>0$ depending only on $M$ 
\[
\left\| \bp \zeta \\ V \ep - \bp g_{+}+g_{-} \\ g_{+}-g_{-} \ep \right\|_{L_{t}^{\infty}(0,T;H_{x}^{1}(\R))} \leq (1 + \tfrac{\mu}{\epsilon}) (1+T) (\epsilon + \mu + \sqrt{\epsilon} \sqrt{T} )  C
\]
so that in the case $\mu \lesssim \epsilon$ we get a convergence result when $\epsilon \to 0$.
\end{comment}
A proof of such a result can be adapted from \cite{bona_colin_lannes} or \cite[Section 7.1]{Lannes_ww}.
\er

\subsection{Local existence}\label{localexist1d}

In this section we provide a local wellposedness result of \eqref{boussinesq_eq1d} on an existence time independent of $\epsilon,\mu \in (0,1]$. We introduce the functional space
\[
X^{3}_{\mu}(\R) := \left\{ (\zeta,V) \in (H^{3} \times H^{4})(\R) \text{  ,  } \| (\zeta,V) \|_{X_{\mu}^{3}} := \| \zeta \|_{H^{3}} + \| V \|_{H^{3}} + \sqrt{\mu} \| \partial_{x} V \|_{H^{3}} < \infty \right\}.
\]

\bpr\label{localexistence1d}
Let $h_{0}>0$, $A>0$, $\epsilon \in (0,1]$ and $\mu \in (0,1]$. Let $(\zeta_{0},V_{0}) \in X^{3}_{\mu}(\R)$ with $1+\epsilon \zeta_{0} \geq h_{0}$ and $\| (\zeta_{0}, V_{0}) \|_{X^{3}_{\mu}} \leq A$. There exists a time $T>0$ that only depends on $A$ and $h_{0}$ and a unique solution $(\zeta,V) \in \mathcal{C}([0,T];X^{3}_{\mu}(\R))$ of \eqref{boussinesq_eq1d} with initial datum $(\zeta_{0},V_{0})$. Furthermore there exists a constant $M>0$ depending only on $A$ and $h_{0}$ such that
\[
\| (\zeta,V) \|_{L^{\infty}(0,T;X^{3}_{\mu})} \leq M.
\]
\epr
The proof is similar to the proofs of \cite[Proposition 2.15]{my_long_wave_corio} or \cite[Theorem 1]{israwi_green_naghdi}. The key point is that System \eqref{boussinesq_eq1d} under the assumption that $(1+\epsilon \zeta) \geq \frac{h_{0}}{2}$ has a symmetrizer
\begin{equation*}
S(U) = \begin{pmatrix} 1 & 0 \\ 0 & (1+\epsilon \zeta) - \frac{\mu}{3} \partial_{x} \left( (1+\epsilon \zeta) \partial_{x} \cdot \right) \end{pmatrix}
\end{equation*}
and corresponding energies 
\begin{equation*}
\mathcal{E}^{k}(U) = \left(S(U) \partial_{x}^{k} U, \partial_{x}^{k} U \right)_{2}
\end{equation*}
that allow to control the $X^{3}_{\mu}(\R)$-norm.

\subsection{Dispersive estimates}

The phase $g:=r \mapsto \frac{r}{\sqrt{1+ \frac{r^2}{3}}}$ is smooth and
\colr
\be\label{phasis_prop}
g'>0 \text{, } g''<0 \text{ on } \R^{+}_{\ast} \text{, } g'(r) -1 \underset{r \sim 0}{\sim}  - \frac{r^2}{2}  \text{, } g''(r) \underset{r \sim 0}{\sim} -r \text{, } g'(r) \underset{r \sim + \infty}{\sim} 3^{\frac32} r^{-3} \text{, } g''(r) \underset{r \sim + \infty}{\sim} 3^{\frac52} r^{-4}.
\ee
\colb
Thanks to Appendix \ref{s:dispestim} one can prove several dispersive estimates. In the following $\chi$ is a smooth compactly supported function that is equal to $1$ near $0$.

Firstly using Lemma \ref{decay1d}\footnote{\colr with $\beta=1,s=\frac12,\alpha=-4,l=2$.\colb}, there exists a constant $C>0$ such that for any Schwartz class function $f$, any $t \neq 0$, any $\mu \in (0,1]$ and any $\epsilon \in (0,1]$
\begin{align*}
&\left\| e^{\pm \frac{t}{\epsilon} \tfrac{\partial_{x}}{\sqrt{1 + \frac{\mu}{3} D^2}}} \chi( \sqrt{\mu} |D|) \sgn(D) |D|^{\frac12} f \right\|_{L^{\infty}_{x}} \leq \frac{C}{\mu^{\frac12}} \left( \frac{\epsilon}{|t|} \right)^{\frac12} \left\| \chi( \sqrt{\mu} |D|) f \right\|_{L^{1}},\\
&\left\| e^{\pm \frac{t}{\epsilon} \tfrac{\partial_{x}}{\sqrt{1 + \frac{\mu}{3} D^2}}} (1-\chi( \sqrt{\mu} |D|)) f \right\|_{L^{\infty}_{x}} \leq \frac{C}{\mu^{\frac12}} \left( \frac{\epsilon}{|t|} \right)^{\frac12} \| (1-\chi( \sqrt{\mu} |D|)) (\sqrt{\mu} |D|) ^2 f \|_{L^{1}}.
\end{align*}
Then corresponding Strichartz estimates can be obtained: for any $(q,r) \in \{ (4,\infty) , (\infty,2) \}$ and any functions $f$ and $G$ smooth enough
\bes
\begin{aligned}
\left\| e^{\pm \frac{t}{\epsilon} \tfrac{\partial_{x}}{\sqrt{1 + \frac{\mu}{3} D^2}}} \partial_{x} f \right\|_{L^{q}_{t} L^{r}_{x}} \lesssim & \left( \frac{\epsilon}{\mu} \right) ^{\frac12 (\frac12 - \frac{1}{r}) } \| |D|^{\frac34 + \frac{1}{2r}} \chi( \sqrt{\mu} |D|) f \|_{L^{2}}\\
&+ \left( \frac{\epsilon}{\mu} \right) ^{\frac12 (\frac12 - \frac{1}{r}) } \| (1-\chi( \sqrt{\mu} |D|)) (\sqrt{\mu} |D|)^{2 (\frac12 - \frac{1}{r})} \partial_{x} f \|_{L^{2}}
\end{aligned}
\ees
\be\label{retardedstrichartz1/21d}
\begin{aligned}
\left\| \int_{0}^{t} e^{\pm \frac{(t-s)}{\epsilon} \tfrac{\partial_{x}}{\sqrt{1 + \frac{\mu}{3} D^2}}} \partial_{x} G(s) ds  \right\|_{L^{q}_{t} L^{r}_{x}} \lesssim & \left( \frac{\epsilon}{\mu} \right)^{\frac12 (1 - \frac{1}{r}) } \| |D|^{\frac12 + \frac{1}{2r}} \chi( \sqrt{\mu} |D|) G \|_{L^{\frac43}_{t} L^{1}_{x}}\\
&+ \left( \frac{\epsilon}{\mu} \right)^{\frac12 (1 - \frac{1}{r}) } \| (1-\chi( \sqrt{\mu} |D|)) (\sqrt{\mu} |D|)^{2 (1 - \frac{1}{r})} \partial_{x} G \|_{L^{\frac43}_{t} L^{1}_{x}}.
\end{aligned}
\ee
\colr One can prove such bounds by applying a $T^{\ast} T$ argument\footnote{\colr in the spirit of \cite[Lemma 2.1]{GV92}.\colb} on the following operators that are both defined from $L^{q'}_{t}(\R,L^{r'}_{x}(\R))$ to $L^{2}(\R)$
\begin{align*}
&T_{LF}(F) = \int_{\R} e^{\mp \frac{s}{\epsilon} \tfrac{\partial_{x}}{\sqrt{1 + \frac{\mu}{3} D^2}}} \chi(\sqrt{\mu} |D|) \frac{\partial_{x}}{|D|^{\frac34 + \frac{1}{2r}}} F(s) ds,\\
&T_{HF}(F) = \int_{\R} e^{\mp \frac{s}{\epsilon} \tfrac{\partial_{x}}{\sqrt{1 + \frac{\mu}{3} D^2}}} \frac{1-\chi(\sqrt{\mu} |D|)}{(\sqrt{\mu} |D|)^{1-\frac{2}{r}}} F(s) ds.
\end{align*}
In a similar way \colb one can also get from Lemma \ref{decay1d}\footnote{\colr with $\beta=1,s=0,\alpha=-4,l=3$.\colb} together with Lemma \ref{control_Bessel} that for any functions $f$ and $F$ smooth enough
\[
\left\| e^{\pm \frac{t}{\epsilon} \tfrac{\partial_{x}}{\sqrt{1 + \frac{\mu}{3} D^2}}} f \right\|_{L^{\infty}} \lesssim \frac{C}{\mu^{\frac13}} \left( \frac{\epsilon}{|t|} \right)^{\frac13} \| (1+\mu D^{2})^{\frac56} f \|_{L^{1}}
\]
and corresponding Strichartz estimates can be obtained from a $T^{\ast} T$ argument: for any $(\tilde q, \tilde r) \in \{ (6,\infty) , (\infty,2) \}$
\be\label{strichartz1/31d}
\left\| e^{\pm \frac{t}{\epsilon} \tfrac{\partial_{x}}{\sqrt{1 + \frac{\mu}{3} D^2}}} f  \right\|_{L^{\tilde q}_{t} L^{\tilde r}_{x}} \lesssim \left(  \frac{\epsilon}{\mu} \right)^{\frac13 (\frac12 - \frac{1}{\tilde r})} \| (1+\mu D^{2})^{\frac56 (\frac12 - \frac{1}{\tilde r})} f \|_{L^{2}}
\ee
\be\label{retardedstrichartz1/31d}
\begin{aligned}
\left\| \int_{0}^{t} e^{\pm \frac{(t-s)}{\epsilon} \tfrac{\partial_{x}}{\sqrt{1 + \frac{\mu}{3} D^2}}} F(s) ds  \right\|_{L^{\tilde q}_{t} L^{\tilde r}_{x}} \lesssim \left(  \frac{\epsilon}{\mu} \right)^{\frac13 (1 - \frac{1}{\tilde r})}  \| (1+\mu D^{2})^{\frac56 (1 - \frac{1}{\tilde r})} F \|_{L^{\frac65}_{t} L^{1}_{x}}.
\end{aligned}
\ee

\subsection{Proof of Theorem \ref{rigidlid1d}}
From Duhamel's principle on \eqref{boussinesq_eq1dcompactform}
\bes
\bp \zeta \\ V \ep(t) - e^{-\frac{t}{\epsilon} A(\partial_{x})} \bp \zeta_{0} \\ V_{0} \ep = \int_{0}^{t} \exp(\tfrac{\tau-t}{\epsilon} A(\partial_{x})) F \left( \bp \zeta \\ V \ep(\tau) \right) d\tau.
\ees
Then from \eqref{exp(tA)1d} and \colr since $F$ is a derivative \colb we can use \eqref{retardedstrichartz1/21d} so that for any $(q,r) \in \{ (4,\infty) , (\infty,2) \}$
\[
\left\|  \int_{0}^{t} \exp(\tfrac{\tau-t}{\epsilon} A(\partial_{x})) F \left( \bp \zeta \\ V \ep(\tau) \right) d\tau \right\|_{L_{t}^{q}(0,T;L_{x}^{r}(\R))} \lesssim  \left( \frac{\epsilon}{\mu} \right)^{\frac14 + \frac{1}{q}} \left( A_{LF} + A_{HF} \right)
\]
where
\[
\begin{aligned}
A_{LF} = &\left\| |D|^{\frac12+\frac{1}{2r}} \left( \zeta V \right) \right\|_{L^{\frac43}_{t}(0,T;L^{1}_{x})} + \left\| \frac{|D|^{\frac12+\frac{1}{2r}}}{\sqrt{1+ \frac{\mu}{3} |D|^2}}\left( V^{2} \right) \right\|_{L^{\frac43}_{t}(0,T;L^{1}_{x})}\\
&+ \left\| \frac{|D|^{\frac12+\frac{1}{2r}} }{\sqrt{1+ \frac{\mu}{3} |D|^2}} \left( \zeta V \right) \right\|_{L^{\frac43}_{t}(0,T;L^{1}_{x})} + \left\| \frac{|D|^{\frac12+\frac{1}{2r}} }{1+ \frac{\mu}{3} |D|^2} \left( V^{2} \right) \right\|_{L^{\frac43}_{t}(0,T;L^{1}_{x})}
\end{aligned}
\]
and
\[
\begin{aligned}
A_{HF} = &\left\| (\sqrt{\mu} |D|)^{2(1-\frac1r)} \partial_{x} \left( \zeta V \right) \right\|_{L^{\frac43}_{t}(0,T;L^{1}_{x})} + \left\| \frac{(\sqrt{\mu} |D|)^{2(1-\frac1r)}}{\sqrt{1+ \frac{\mu}{3} |D|^2}} \partial_{x} \left( V^{2} \right) \right\|_{L^{\frac43}_{t}(0,T;L^{1}_{x})}\\
&+ \left\| \frac{(\sqrt{\mu} |D|)^{2(1-\frac1r)}}{\sqrt{1+ \frac{\mu}{3} |D|^2}} \partial_{x} \left( \zeta V \right) \right\|_{L^{\frac43}_{t}(0,T;L^{1}_{x})} + \left\| \frac{(\sqrt{\mu} |D|)^{2(1-\frac1r)}}{1+ \frac{\mu}{3} |D|^2} \partial_{x} \left( V^{2} \right) \right\|_{L^{\frac43}_{t}(0,T;L^{1}_{x})}.
\end{aligned}
\]
Lemmas \ref{control_Bessel} and \ref{control|D|1d} provide
\[
A_{LF} + A_{HF} \lesssim \| \| \zeta \|_{H^{3}_{x}} \| V \|_{H^{3}_{x}} \|_{L^{\frac43}_{t}(0,T)} + \| \| V \|_{H^{3}_{x}}^2 \|_{L^{\frac43}_{t}(0,T)} \lesssim T^{\frac34} M^2.
\]
We then get the first bound when $(q,r) \in \{ (4,\infty) , (\infty,2) \}$ and the other cases follow from H\"older's inequality. The second bound follows the same way using instead \eqref{strichartz1/31d} and \eqref{retardedstrichartz1/31d}. \colr Note that since the initial datum is not necessary a derivative one can not use the homogeneous version of \eqref{retardedstrichartz1/21d}, hence a weaker rate of converge in that case\colb. Finally the third bound follows from Morawetz-type estimates established in Proposition \ref{morawetz}.

\section{The classical 2D Boussinesq system}\label{s:classicalBoussi2d}

The two dimensional classical Boussinesq system reads as

\be\label{boussinesq_eq2d}
\left\{
\begin{aligned}
&\epsilon \partial_{t} \zeta + \nabla\cdot \left( \left[1 + \epsilon \zeta \right] \bV \right) = 0,\\
&\epsilon \left(1- \frac{\mu}{3} \nabla \nabla \cdot \right) \partial_{t} \bV + \nabla \zeta + \epsilon (\bV \cdot \nabla) \bV = 0.
\end{aligned}
\right.
\ee
Taking the divergence of the second equation and denoting by $\bU := (\zeta,\nabla \cdot \bV)^{T} \in \R^{2}$, we get the following system
\bes
\epsilon \partial_{t} \bU + A(D) \bU = \epsilon F(\zeta,\bV)
\ees
where
\bes
A(D) = \begin{pmatrix} 0 & 1 \\ (1- \frac{\mu}{3} \Delta)^{-1} \Delta & 0 \end{pmatrix} \text{   and   } F(\zeta,\bV) = - \begin{pmatrix} \nabla \cdot (\zeta \bV) \\ (1- \frac{\mu}{3} \Delta)^{-1} \nabla \cdot ( (\bV \cdot \nabla ) \bV) \end{pmatrix}. 
\ees
Note that
\bes
\exp(tA(D)) = \begin{pmatrix} \cos \left( t \frac{|D|}{\sqrt{1+ \frac{\mu}{3} |D|^2}} \right) & \frac{\sqrt{1+ \frac{\mu}{3} |D|^2}}{|D|} \sin \left( t \frac{|D|}{\sqrt{1+ \frac{\mu}{3} |D|^2}} \right) \\ -\frac{|D|}{\sqrt{1+ \frac{\mu}{3} |D|^2}} \sin \left( t \frac{|D|}{\sqrt{1+ \frac{\mu}{3} |D|^2}} \right) & \cos \left( t \frac{|D|}{\sqrt{1+ \frac{\mu}{3} |D|^2}} \right) \end{pmatrix}.
\ees
As in the 1d case we expect $(\zeta,\nabla \cdot \bV)$ to weakly converge to $0$. We however have no control on the rotational component of $\bV$. Applying the operator $\nabla^{\perp} \cdot$ to the second equation and denoting by $\omega := \nabla^{\perp} \cdot \bV$ we get the following equation
\[
\partial_{t} \omega + (\bV \cdot \nabla) \omega + (\nabla \cdot \bV) \omega = 0.
\]
Since $\bV$ is a vector field on $\R^{2}$ it has a Hodge-Weyl decomposition 
\[
\bV = \nabla \frac{\nabla}{\Delta} \cdot \bV + \nabla ^{\perp}\frac{\nabla^{\perp}}{\Delta} \cdot \bV. 
\]
Therefore we expect $\omega$ to converge to $\tilde \omega$ as $\epsilon \to 0$ where 
\[
\partial_{t} \tilde \omega + \left( \nabla ^{\perp}\frac{\nabla^{\perp}}{\Delta} \tilde \omega \cdot \nabla \right) \tilde \omega = 0.
\]
The previous equation is the vorticity formulation of the incompressible Euler equation
\be\label{incompEuler}
\partial_{t} \tilde \bV + (\tilde \bV \cdot \nabla) \tilde \bV + \nabla P = 0 \text{  ,  } \nabla \cdot \tilde \bV = 0
\ee
with $\tilde \omega = \nabla^{\perp} \cdot \tilde \bV$. One can now state the main results of this section.

\bt\label{rigidlid2d1}
Let $M>0$, $T>0$, $\epsilon \in (0,1]$ and $\mu \in (0,1]$. Let $(\zeta,\bV) \in \mathcal{C}([0,T];(H^{6} \times H^{6})(\R^2))$ be a solution of \eqref{boussinesq_eq2d} with initial datum $(\zeta_{0},\bV_{0})$ and  $\tilde \bV \in \mathcal{C}([0,T];L^{2}(\R^2))$ be a solution of the incompressible Euler equation \eqref{incompEuler} with initial datum $\nabla^{\perp} \frac{\nabla^{\perp}}{\Delta} \cdot \bV_{0}$ such that
\[
\| (\zeta,\bV) \|_{L^{\infty}(0,T;(H^{6} \times H^{6})(\R^2))} + \| \tilde \bV \|_{L^{\infty}(0,T;L^{2}(\R^2))} \leq M.
\]
There exists a constant $C>0$ \colr independent of $M$, $T$, $\epsilon$ and $\mu$ \colb such that for any $q,r \geq 2$ with $\frac{1}{q} + \frac{1}{r} = \frac{1}{2}$
\[
\begin{aligned}
&\left\| \bp \zeta \\ \frac{\nabla \nabla}{\Delta} \cdot \bV \ep - \bp 1 & 0 \\ 0 & \frac{\partial_{1}}{\Delta} \\ 0 & \frac{\partial_{2}}{\Delta}  \ep \exp(-\tfrac{t}{\epsilon}A(D)) \bp \zeta_{0} \\ \nabla \cdot \bV_{0} \ep \right\|_{L_{t}^{q}(0,T;L_{x}^{r}(\R^2))} \hspace{-0.5cm}\leq \left( \frac{\epsilon}{\sqrt{\mu}}  \ln(\colr 1 \colb + \tfrac{\mu}{\epsilon^2} T) \right)^{\frac12 + \frac{1}{q}} \colr M^2 \colb T^{\frac12} C\\
&\hspace{11cm} + \tfrac{\epsilon}{\sqrt{\mu}} \colr M^2 \colb T C,\\
&\left\| \bp \zeta \\ \frac{\nabla \nabla }{\Delta} \cdot \bV \ep   \right\|_{L_{t}^{q}(0,T;L_{x}^{r}(\R^2))} \leq \left( \frac{\epsilon}{\sqrt{\mu}}  \ln(\colr 1 \colb + \tfrac{\mu}{\epsilon^2} T) \right)^{\frac{1}{q}} \left( \colr M \colb  + \colr M^2 \colb \left( \frac{\epsilon}{\sqrt{\mu}}  \ln( \colr 1 \colb + \tfrac{\mu}{\epsilon^2} T) T \right)^{\frac12} \right) C\\
&\hspace{5cm} +\tfrac{\epsilon}{\sqrt{\mu}} (\colr M \colb + \colr M^2 \colb T ) C,\\
&\colp \left\| \nabla^{\perp} \frac{\nabla^{\perp}}{\Delta}  \cdot \bV - \tilde \bV \right\|_{L_{t}^{\infty}(0,T;L_{x}^{2}(\R^2))} \hspace{-0.6cm} \leq \left( \frac{\epsilon}{\sqrt{\mu}} \ln( 1 + \tfrac{\mu}{\epsilon^2} T) \right)^{\hspace{-0.1cm} \frac12} \hspace{-0.2cm} \left( \hspace{-0.1cm} M \hspace{-0.05cm} + \hspace{-0.05cm} M^2 \hspace{-0.1cm} \left( \hspace{-0.1cm} \frac{\epsilon}{\sqrt{\mu}} \ln( 1 + \tfrac{\mu}{\epsilon^2} T) T \right)^{\hspace{-0.1cm} \frac12} \hspace{-0.1cm} \right) \hspace{-0.1cm} M\sqrt{T} e^{C M  T} \hspace{-0.05cm}  C\\
&\hspace{6cm} \colp+\tfrac{\epsilon}{\sqrt{\mu}} ( M + M^2 T) M\sqrt{T} e^{C M  T}  C \colb.
\end{aligned}
\]
\et 

\begin{comment}
\br 
Unfortunately one can not take $\tilde q = 2$ in the $\zeta$ part. The main issue is that the Riesz transform is not bounded on $L^1$, see Remark \ref{notL1} in the proof of the theorem.
\er
\end{comment}

\br 
$(\zeta,\frac{\nabla}{\Delta} \nabla \cdot \bV )$ weakly converges to $0$ as $\epsilon \to 0$ whereas $\frac{\nabla^{\perp}}{\Delta} \nabla^{\perp} \cdot \bV$ strongly converges to a solution of the incompressible Euler equation. Therefore $\bV$ strongly converges to a solution of the incompressible Euler equation if and only if $(\zeta_{0} , \nabla \cdot \bV_{0})=(0,0)$. The default of strong convergence is exhibited through the first bound of the theorem. Finally unlike the 1d case, if $\epsilon \sim \mu$ as in \cite{bona_chen_saut_derivation,bona_colin_lannes}, the theorem provides a convergence result as $\epsilon \to 0$.
\er

The previous theorem does not provide uniform bounds with respect to $\mu \to 0$. It is the purpose of the following theorem.

\bt\label{rigidlid2d2}
Let $M>0$, $T>0$, $\epsilon \in (0,1]$ and $\mu \in (0,1]$. Let $(\zeta,\bV) \in \mathcal{C}([0,T];(H^{6} \times H^{6})(\R^2))$ be a solution of \eqref{boussinesq_eq2d} with initial datum $(\zeta_{0},\bV_{0})$ and $\tilde \bV \in \mathcal{C}([0,T];L^{2}(\R^2))$ be a solution of the incompressible Euler equation \eqref{incompEuler} with initial datum $\nabla^{\perp} \frac{\nabla^{\perp}}{\Delta} \cdot \bV_{0}$ such that
\[
\| (\zeta,\bV) \|_{L^{\infty}(0,T;(H^{6} \times H^{6})(\R^2))} + \| \tilde \bV \|_{L^{\infty}(0,T;L^{2}(\R^2))} \leq M.
\]
There exists a constant $C>0$ \colr independent of $M$, $T$, $\epsilon$ and $\mu$ \colb such that for any $q,r \geq 2$ with $\frac{1}{q} + \frac{1}{2r} = \frac{1}{4}$
\[
\begin{aligned}
&\left\| \bp \zeta \\ \frac{\nabla \nabla}{\Delta} \cdot \bV \ep - \bp 1 & 0 \\ 0 & \frac{\partial_{1}}{\Delta} \\ 0 & \frac{\partial_{2}}{\Delta}  \ep \exp(-\tfrac{t}{\epsilon}A(D)) \bp \zeta_{0} \\ \nabla \cdot \bV_{0} \ep \right\|_{L_{t}^{q}(0,T;L_{x}^{r}(\R^2))} \leq \epsilon^{\frac14 + \frac{1}{q}} \colr M^2 \colb T^{\frac34} C,\\
&\left\| \bp \zeta \\ \frac{\nabla \nabla }{\Delta} \cdot \bV \ep \right\|_{L_{t}^{q}(0,T;L_{x}^{r}(\R^2))} \leq \epsilon^{\frac{1}{q}} \left( \colr M \colb  + \colr M^2 \colb T^{\frac34} \epsilon^{\frac14} \right) C,\\
&\colp\left\| \nabla^{\perp} \frac{\nabla^{\perp}}{\Delta}  \cdot \bV - \tilde \bV \right\|_{L_{t}^{\infty}(0,T;L_{x}^{2}(\R^2))} \leq \epsilon^{\frac{1}{4}} \left( M  + M^2  T^{\frac34} \epsilon^{\frac14}  \right) M T^{\frac34} e^{C M T} C,\colb\\
&\sup_{x_{0} \in \R^2} \left\| e^{-(x-x_{0})^2} \bp \zeta \\ \nabla \frac{\nabla}{\Delta} \cdot \bV \ep \right\|_{L_{t}^{2}(0,T;L^{2}_{x}(\R^2))} \leq \epsilon^{\frac12}(\colr M \colb + \colr M^2 \colb  T) C.
\end{aligned}
\]
\et 

\br 
The last two estimates provide the convergence  in $L^{2}_{t}(0,T;L^{2}_{loc}(\R^2))$. We then get  the convergence of the local energy to the local energy of the corresponding solution of the incompressible Euler equation. If furthermore $\tilde \bV \in \mathcal{C}([0,T];H^{6}(\R^2))$ then the convergence is in $L^{2}_{t}(0,T;H^{s}_{loc}(\R^2))$ for any $s \in [0,6)$. 
\er

\subsection{Local existence}\label{localexist2d}

In this section we provide a local wellposedness result of \eqref{boussinesq_eq2d} on an existence time independent of $\epsilon,\mu \in (0,1]$. Let $k \in \N$ with $k \geq 3$. We introduce the functional space
\[
X^{k}_{\mu}(\R^{2}) := \left\{ (\zeta,V) \in (H^{k} \times H^{k})(\R^{2}) \text{  ,  } \| (\zeta,V) \|_{X_{\mu}^{k}} := \| \zeta \|_{H^{k}} + \| V \|_{H^{k}} + \sqrt{\mu} \| \nabla \cdot V \|_{H^{k}} < \infty \right\}.
\]
\bpr\label{localexistence2d}
Let $k \in \N$, $k \geq 3$. Let $h_{0}>0$, $A>0$, $\epsilon \in (0,1]$ and $\mu \in (0,1]$. Let $(\zeta_{0},V_{0}) \in X^{k}_{\mu}(\R^2)$ with $1+\epsilon \zeta_{0} \geq h_{0}$ and $\| (\zeta_{0}, V_{0}) \|_{X^{k}_{\mu}} \leq A$. There exists a time $T>0$ that only depends on $A$ and $h_{0}$ and a unique solution $(\zeta,V) \in \mathcal{C}([0,T];X^{k}_{\mu}(\R^2))$ of \eqref{boussinesq_eq2d} with initial datum $(\zeta_{0},V_{0})$. Furthermore there exists a constant $M>0$ depending only on $A$ and $h_{0}$ such that
\[
\| (\zeta,V) \|_{L^{\infty}(0,T;X_{\mu}^{k})} \leq M.
\]
\epr
The proof is similar to the proofs of \cite[Proposition 6.7]{Lannes_ww} or \cite[Proposition 2.7]{my_KPcorio}. Again the key point is that System \eqref{boussinesq_eq2d} under the assumption that $(1+\epsilon \zeta) \geq \frac{h_{0}}{2}$ has a symmetrizer and corresponding energies that allow to control the $X^{k}_{\mu}(\R^2)$-norm.

\subsection{Dispersive estimates}

From the properties \eqref{phasis_prop} together with the fact $|g'''(r)| \lesssim r^{-5}$, and thanks to Appendix \ref{s:dispestim} one can prove several dispersive estimates.

Firstly, using Lemma \ref{decay2d}\footnote{\colr with $\beta=1$, $s=\alpha=-4$. \colb} together with Lemma \ref{control_Bessel}, for any $m \in \N$ there exists a constant $C>0$ such that for any Schwartz class function $f$, any $t \neq 0$, any $\mu \in (0,1]$ and any $\epsilon \in (0,1]$ we have 
\bes
\left\| e^{\pm \iD \frac{t}{\epsilon} \tfrac{|D|}{\sqrt{1 + \frac{\mu}{3} D^2}}} \frac{\nabla^m}{|D|^m} f \right\|_{L^{\infty}_{x}} \leq \frac{C}{\sqrt{\mu}} \frac{\epsilon}{|t|} \| (1+\mu |D|^2)^{2} f \|_{L^{1}}.
\ees
Then corresponding Strichartz estimates can be obtained from a $T^{\ast} T$ argument: for any $q,r,\tilde q,\tilde r \geq 2$ such that $\frac{1}{q} + \frac{1}{r} = \frac12$,  $\frac{1}{\tilde q} + \frac{1}{\tilde r} = \frac12$  with $r,\tilde r<\infty$, any $m \in \N$, any $\epsilon \in (0,1]$ and any $\mu \in (0,1]$
\begin{equation}\label{strichartz2d}
\begin{aligned}
&\left\| e^{\pm \iD \frac{t}{\epsilon} \tfrac{|D|}{\sqrt{1 + \frac{\mu}{3} D^2}}} \frac{\nabla^{m}}{|D|^{m}} f  \right\|_{L^{q}_{t} L^{r}_{x}} \lesssim \left( \frac{\epsilon}{\sqrt{\mu}} \right)^{\frac12 - \frac{1}{r}} \| (1+\mu |D|^2)^{2 (\frac12 - \frac{1}{r})} f \|_{L^{2}}\\
&\left\| \int_{0}^{t} e^{\pm \iD \frac{(t-s)}{\epsilon} \tfrac{|D|}{\sqrt{1 + \frac{\mu}{3} D^2}}} \frac{\nabla^{m}}{|D|^{m}} F(s) ds  \right\|_{L^{q}_{t} L^{r}_{x}} \lesssim \left( \frac{\epsilon}{\sqrt{\mu}} \right)^{1 - \frac{1}{\tilde r} - \frac{1}{r}} \| (1+\mu |D|^2)^{2(1 - \frac{1}{\tilde r} - \frac{1}{r})} F \|_{L^{\tilde q'}_{t} L^{\tilde r'}_{x}}.
\end{aligned}
\end{equation}
It is well-known that the previous estimates do not work at the endpoints $(q,r)=(2,\infty)$ or $(\tilde q, \tilde r)=(2,\infty)$. We can however prove a logarithmic estimate for functions whose Fourier transform are well localized. Such type of estimates were performed for the wave equation (\cite{JMR00} or \cite[Theorem 8.30]{BCD_Fourier}) or the Schr\"odinger equation (\cite{Tao_counterexample}). We provide in the following a general result.

\bpr\label{propstrichartz12d}
Let $\chi$ be a smooth compactly supported function, $g$ a function defined on $\R^{\ast}_{+}$, $\alpha \in \R$ and $m \in \N$. Assume there exists a constant $C_{0}>0$ such that for any Schwartz class function $f$ and any $t \neq 0$
\bes
\left\| e^{\pm \iD t g(|D|)} \frac{\nabla^{m}}{|D|^{m}} f \right\|_{L^{\infty}_{x}} \leq \frac{C_{0}}{|t|} \| (1+|D|^2)^{\frac{\alpha}{2}} f \|_{L^{1}}.
\ees
Then there exists a constant $C>0$ such that for any $\lambda>0$, any $T>0$, any $\mu>0$ and any $\epsilon>0$ we have
\[
\left\| e^{\pm \iD \frac{t}{\epsilon \sqrt{\mu}} g(\sqrt{\mu} |D|)}\chi(\lambda^{-1} |D|)  \frac{\nabla^{m}}{|D|^{m}} f \right\|_{L^{2}_{t}(0,T;L^{\infty}_{x})} \leq  C \left( \frac{\epsilon}{\sqrt{\mu}} \right)^{\frac12} \sqrt{\ln \left( 1 + \frac{\sqrt{\mu}}{\epsilon} \lambda^2 T \right) } \| (1+\mu |D|^{2})^{\frac{\alpha}{4}} f \|_{L^{2}}
\]
and if we define the operator $H$ as
\[
H(F) := \int_{0}^{t} e^{\pm \iD \frac{(t-s)}{\epsilon \sqrt{\mu}} g(\sqrt{\mu} |D|)} \chi(\lambda^{-1} |D|)  \frac{\nabla^{m}}{|D|^{m}} F(s,\cdot) ds
\]
we have
\[
\begin{aligned}
&\left\| H(F) \right\|_{L^{\infty}_{t}(0,T;L^{2}_{x})} \leq C \left( \frac{\epsilon}{\sqrt{\mu}} \right)^{\frac12} \sqrt{\ln \left( 1 + \frac{\sqrt{\mu}}{\epsilon} \lambda^2 T \right) } \| (1+\mu |D|^{2})^{\frac{\alpha}{4}} F \|_{L^{2}_{t}(0,T;L^{1}_{x})},\\
&\left\|  H(F) \right\|_{L^{2}_{t}(0,T;L^{\infty}_{x})} \leq  C \left( \frac{\epsilon}{\sqrt{\mu}} \right)^{\frac12}  \sqrt{\ln \left( 1 + \frac{\sqrt{\mu}}{\epsilon} \lambda^2 T \right) } \| (1+\mu |D|^{2})^{\frac{\alpha}{4}} F \|_{L^{1}_{t}(0,T;L^{2}_{x})},\\
&\left\|  H(F)  \right\|_{L^{2}_{t}(0,T;L^{\infty}_{x})} \leq  C \frac{\epsilon}{\sqrt{\mu}}  \ln \left( 1 + \frac{\sqrt{\mu}}{\epsilon} \lambda^2 T \right) \| (1+\mu |D|^{2})^{\frac{\alpha}{2}} F \|_{L^{2}_{t}(0,T;L^{1}_{x})}.
\end{aligned}
\]
\epr

\begin{proof}
By assumption and change of variables we actually have for any $\mu \in (0,1]$ and any $\epsilon \in (0,1]$
\be\label{assumppropstrichartz12d}
\left\| e^{\pm \iD \frac{t}{\epsilon \sqrt{\mu}} g(\sqrt{\mu} |D|)} \frac{\nabla^{m}}{|D|^{m}} f \right\|_{L^{\infty}_{x}} \leq \frac{\epsilon}{\sqrt{\mu}} \frac{C_{0}}{|t|} \| (1+\mu |D|^2)^{\frac{\alpha}{2}} f \|_{L^{1}}.
\ee
We only prove the case $m=0$ since the methodology is the same for the other cases. We introduce the operator 
\[
\begin{array}{ccccc}
S & : & L^{2}_{t}(0,T, L^{1}_{x}(\R^{n})) & \to & L^{2}(\R^{2})\\
 & & F & \mapsto & \displaystyle{\int_{0}^{T}} e^{\mp \iD \frac{s}{\epsilon \sqrt{\mu}} g(\sqrt{\mu} |D|)} (1+\mu |D|^2)^{-\frac{\alpha}{4}} \chi(\lambda^{-1} |D|)  F(s) ds.
\end{array}
\]
For any $f \in L^{2}(\R^{2})$
\begin{align*}
\| S^{\ast} f \|_{L^{2}_{t}(0,T;L^{\infty}_{x})} &\leq  \sup_{\| G \|_{L^{2}_{t}L^{1}_{x}} \leq 1} \left\{ \int_{0}^{T} \left< e^{\mp \iD \frac{t}{\epsilon \sqrt{\mu}} g(\sqrt{\mu} |D|)} (1+\mu |D|^2)^{-\frac{\alpha}{4}} \chi(\lambda^{-1} |D|)  f , G(t, \cdot ) \right>_{L^{2}_{x} \times L^{2}_{x}} dt \right\}\\
&\leq \| f \|_{L^{2}} \sup_{\| G \|_{L^{2}_{t}(0,T;L^{1}_{x})} \leq 1} J_{G}
\end{align*}
where
\[
J_{G} := \left\{ \left\| \int_{0}^{T} e^{\mp \iD \frac{t}{\epsilon \sqrt{\mu}} g(\sqrt{\mu} |D|)} (1+\mu |D|^2)^{-\frac{\alpha}{4}} \chi(\lambda^{-1} |D|)  G(t,\cdot) dt \right\|_{L^{2}_{x}}  \right\}.
\]
Then
\[
J_{G}^2 \leq \int_{[0,T]^2} \underbrace{ \left\| e^{\mp \iD \frac{(t-s)}{\epsilon \sqrt{\mu}} g(\sqrt{\mu} |D|)} (1+\mu |D|^2)^{-\frac{\alpha}{2}} \chi^2(\lambda^{-1} |D|)  G(t,\cdot) \right\|_{L^{\infty}_{x}} }_{:=K_{G}}  \| G(s,\cdot) \|_{L^{1}_{x}} dt ds
\]
and using \eqref{assumppropstrichartz12d}
\[
K_{G} \lesssim \frac{\epsilon}{\sqrt{\mu}} \frac {1}{|t-s|} \| G(t,\cdot) \|_{L^{1}_{x}}
\]
whereas from Bernstein's Lemma \ref{Bernstein}
\[
\begin{aligned}
K_{G} &\lesssim \lambda \left\| e^{\mp \iD \frac{(t-s)}{\epsilon \sqrt{\mu}} g(\sqrt{\mu} |D|)} (1+ \mu |D|^2)^{-\frac{\alpha}{2}}  \chi(\lambda^{-1} |D|)  G(t,\cdot) \right\|_{L^{2}_{x}}\\
&\lesssim \lambda \left\| \chi^2(\lambda^{-1} |D|)  G(t,\cdot) \right\|_{L^{2}_{x}} \lesssim \lambda^{2} \left\| G(t,\cdot) \right\|_{L^{1}_{x}}
\end{aligned}
\]
so that
\[
J_{G}^2 \leq \int_{[0,T]^2} \frac{\lambda^{2}}{1+ \lambda^{2} \frac{\sqrt{\mu}}{\epsilon} |t-s|} \| G(t,\cdot) \|_{L^{1}_{x}} \| G(s,\cdot) \|_{L^{1}_{x}} dt ds.
\]
The first bound follows from Schur's test. Then for any $(a,b),(\tilde a,\tilde b) \in \{ (2,\infty),(\infty,2) \}$, we define the operator $R_{\tilde a,a} $
\[
\begin{array}{ccccc}
R_{\tilde a,a} & : & L^{\tilde b'}_{t}(\R, L^{\tilde a'}_{x}(\R^{2})) & \to & L^{b}_{t}(\R, L^{a}_{x}(\R^{2}))\\
 & & F & \mapsto & \displaystyle \int_{0}^{t} e^{\pm \iD \frac{(t-s)}{\epsilon \sqrt{\mu}} g(\sqrt{\mu} |D|)} (1+\mu |D|^2)^{-\frac{\alpha}{4}} \chi(\lambda^{-1} |D|) F(s,\cdot) ds.
\end{array}
\]
We note that 
\[
\left\| R_{\infty,2}(F) \right\|_{L^{\infty}_{t}(0,T;L^{2}_{x})} = \sup_{t \in [0,T]} \left\| S ( \mathds{1}_{(0, t)} F) \right\|_{L^{2}_{x}}
\]
so that the second bound follows by duality and the first bound. Furthermore,  we notice that $R_{2,\infty} +R_{\infty,2}^{\ast} = S^{\ast} \tilde S$ where
\[
\begin{array}{ccccc}
\tilde S & : & L^{1}_{t}(0,T, L^{2}_{x}(\R^{n})) & \to & L^{2}(\R^{2})\\
 & & F & \mapsto & \displaystyle{\int_{0}^{T}} e^{\mp \iD \frac{s}{\epsilon \sqrt{\mu}} g(\sqrt{\mu} |D|)}  F(s) ds.
\end{array}
\]
The third estimate follows from the first and second estimates together with the fact $\tilde S$ is bounded since $\tilde S^{\ast}$ is bounded. Finally denoting
\[
L_{F} := \left\| \int_{0}^{t} e^{\pm \iD \frac{(t-s)}{\epsilon \sqrt{\mu}} g(\sqrt{\mu} |D|)} (1+ \mu |D|^2)^{-\frac{\alpha}{2}} \chi(\lambda^{-1} |D|)  F(s,\cdot) ds \right\|_{L^{2}_{t}(0,T;L^{\infty}_{x})} 
\]
and proceeding similarly as for the bound on $K_{G}$ we get
\begin{align*}
L_{F} & \lesssim \left\| \int_{s=0}^{t} \frac{\lambda^{2}}{1+\lambda^{2} \frac{\sqrt{\mu}}{\epsilon} |t-s|} \| F(s,\cdot) \|_{L^{1}_{x}} ds \right\|_{L^{2}_{t}(0,T)}\\
&\lesssim \sqrt{ \sup_{t \in [0,T] } \int_{s=0}^{t} \frac{\lambda^{2}}{1+\lambda^{2} \frac{\sqrt{\mu}}{\epsilon} |t-s|} ds} \sqrt{ \int_{t=0}^{T} \int_{s=0}^{t} \frac{\lambda^{2}}{1+\lambda^{2} \frac{\sqrt{\mu}}{\epsilon} |t-s|} \| F(s,\cdot) \|_{L^{1}_{x}}^2 ds dt }
\end{align*}
and the fourth bound follows.
\end{proof}

Secondly, we get from Lemma \ref{decay2d1/2}(i)\footnote{\colr with $l=2$, $\eta=\frac12$.\colb} together with Lemma \ref{control_Bessel} that if $\chi$ is a smooth compactly supported that is equal to $1$ near $0$ and $m \in \N$
\[
\left\| e^{\pm \iD \frac{t}{\epsilon} \tfrac{|D|}{\sqrt{1 + \frac{\mu}{3} D^2}}} \chi( \sqrt{\mu} |D|) \frac{\nabla^{m}}{|D|^{m}} f \right\|_{L^{\infty}_{x}} \lesssim \frac{\epsilon^{\frac12}}{|t|^{\frac12}}\| (1+|D|^2) f \|_{L^{1}}
\]
\colr and from Lemma \ref{decay2d}(ii)\footnote{\colr with $s=-3$, $\alpha=-4$.\colb} together with Lemma \ref{control_Bessel} that
\[
\left\| e^{\pm \iD \frac{t}{\epsilon} \tfrac{|D|}{\sqrt{1 + \frac{\mu}{3} D^2}}} (1-\chi( \sqrt{\mu} |D|)) \frac{\nabla^{m}}{|D|^{m}} f \right\|_{L^{\infty}_{x}} \lesssim \mu^{\frac34} \frac{\epsilon^{\frac12}}{|t|^{\frac12}}\| (1+|D|^2)^{\frac32} f \|_{L^{1}}.
\]
\colb
One can then obtained corresponding Strichartz estimates from a $T^{\ast} T$ argument: for any $(\tilde q, \tilde r) \in \{ (4,\infty) , (\infty,2) \}$, for any $\mu \in (0,1]$, any $\epsilon \in (0,1]$ and any $m \in \N$
\be\label{strichartz1/22d}
\begin{aligned}
&\left\| e^{\pm \iD \frac{t}{\epsilon} \tfrac{|D|}{\sqrt{1 + \frac{\mu}{3} D^2}}} \frac{\nabla^{m}}{|D|^{m}} f  \right\|_{L^{\tilde q}_{t} L^{\tilde r}_{x}} \lesssim \epsilon^{\frac12(\frac12-\frac{1}{\tilde r})} \| (1+|D|^2)^{\frac32(\frac12-\frac{1}{\tilde r})} f \|_{L^{2}},\\
&\left\| \int_{0}^{t} e^{\pm \iD \frac{t-s}{\epsilon} \tfrac{|D|}{\sqrt{1 + \frac{\mu}{3} D^2}}} \frac{\nabla^{m}}{|D|^{m}} F(s) ds  \right\|_{L^{\tilde q}_{t} L^{\tilde r}_{x}} \lesssim \epsilon^{\frac12(1-\frac{1}{\tilde r})} \| (1+|D|^2)^{\frac32(1-\frac{1}{\tilde r})} F \|_{L^{\frac43}_{t} L^{1}_{x}}.
\end{aligned}
\ee

\subsection{Proofs of Theorem \ref{rigidlid2d1} and Theorem \ref{rigidlid2d2}}
We begin with the proof of Theorem \ref{rigidlid2d1}. Using Duhamel's principle
\bes
\bp \zeta \\ \frac{\nabla}{\Delta} \nabla \cdot \bV \ep - \bp 1 & 0 \\ 0 & \frac{\partial_{1}}{\Delta} \\ 0 & \frac{\partial_{2}}{\Delta}  \ep \exp(-tA(D)) \bp \zeta_{0} \\ \nabla \cdot \bV_{0} \ep = \underbrace{\int_{0}^{t} \bp 1 & 0 \\ 0 & \frac{\partial_{1}}{\Delta} \\ 0 & \frac{\partial_{2}}{\Delta}  \ep \exp(\tfrac{\tau-t}{\epsilon} A(D))  F(\zeta(\tau),\bV(\tau)) d\tau}_{:= I}.
\ees
Let $\frac{1}{q} + \frac{1}{r} = \frac{1}{2}$ and $\chi$ be a smooth compactly supported function that is equal to $1$ near $0$. By Proposition \ref{propstrichartz12d} and interpolation, for any $j \in \{ 1,2 \}$
\[
\left\| \chi(\sqrt{\tfrac{\epsilon}{\sqrt{\mu}}}|D|)  I_{j} \right\|_{L_{t}^{q}(0,T;L_{x}^{r}(\R^2))} \lesssim \left( \frac{\epsilon}{\sqrt{\mu}} \ln(1 + \tfrac{\mu}{\epsilon^2} T) \right)^{\frac{1}{q} + \frac12} B_{j}
\]
where
\[
\begin{aligned}
&B_{1} := \left\| (1+ \mu|D|^2)^{2(1-\frac1r)} \nabla \cdot (\zeta \bV) \right\|_{L^{2}_{t}(0,T;L^{1}_{x})} + \left\| \frac{(1+ \mu|D|^2)^{2(1-\frac1r)}}{\sqrt{1+ \frac{\mu}{3} |D|^2}} ((\bV \cdot \nabla) \bV) \right\|_{L^{2}_{t}(0,T;L^{1}_{x})}\\
&B_{2} := \left\| \frac{(1+ \mu|D|^2)^{2(1-\frac1r)} }{\sqrt{1+ \frac{\mu}{3} |D|^2}} \nabla \cdot \left( \zeta \bV \right) \right\|_{L^{2}_{t}(0,T;L^{1}_{x})} + \left\| \frac{(1+ \mu|D|^2)^{2(1-\frac1r)}}{1+ \frac{\mu}{3} |D|^2} ( (\bV \cdot \nabla) \bV ) \right\|_{L^{2}_{t}(0,T;L^{1}_{x})}.
\end{aligned}
\]
Using Lemmas \ref{control_Bessel} and \ref{control|D|2d} we get
\[
B_{1} + B_{2} \lesssim \| \| \zeta \|_{H^{5}_{x}} \| \bV \|_{H^{5}_{x}} \|_{L^{2}_{t}(0,T)} +\|  \| \bV \|_{H^{5}_{x}}^2 \|_{L^{2}_{t}(0,T)} \lesssim T^{\frac12} M^2.
\]
Furthermore using Lemma \ref{controlHFLp}
\[
\begin{aligned}
\left\| (1-\chi(\sqrt{\tfrac{\epsilon}{\sqrt{\mu}}} |D|))  I_{j} \right\|_{L_{t}^{\infty}(0,T;L_{x}^{2}(\R^2))} &\lesssim \frac{\epsilon}{\sqrt{\mu}} \left\| |D|^{2} I_{j} \right\|_{L_{t}^{\infty}(0,T;L_{x}^{2}(\R^2))}\\
&\lesssim \frac{\epsilon}{\sqrt{\mu}} \left(  \| \| \zeta \|_{H^{3}_{x}} \| \bV \|_{H^{3}_{x}} \|_{L^{1}_{t}(0,T)} +\|  \| \bV \|_{H^{3}_{x}}^2 \|_{L^{1}_{t}(0,T)} \right) \\
&\lesssim T \frac{\epsilon}{\sqrt{\mu}}  M^2
\end{aligned}
\]
whereas from Sobolev inequalities, Lemma \ref{controlHFLp} and \eqref{strichartz2d} with $r=\tilde r = \frac14$
\[
\begin{aligned}
\left\| (1-\chi(\sqrt{\tfrac{\epsilon}{\sqrt{\mu}}} |D|))  I_{j} \right\|_{L_{t}^{2}(0,T;L_{x}^{\infty}(\R^2))} &\lesssim \left\| (1-\chi(\sqrt{\tfrac{\epsilon}{\sqrt{\mu}}} |D|))  I_{j} \right\|_{L_{t}^{2}(0,T;W_{x}^{1,4}(\R^2))}\\
&\lesssim  T^{\frac14} \left\| (1-\chi(\sqrt{\tfrac{\epsilon}{\sqrt{\mu}}} |D|))  I_{j} \right\|_{L_{t}^{4}(0,T;W_{x}^{1,4}(\R^2))}\\
&\lesssim T^{\frac14} \left( \frac{\epsilon}{\sqrt{\mu}} \right)^{\frac12} \left\| I_{j} \right\|_{L_{t}^{4}(0,T;W_{x}^{2,4}(\R^2))}\\
&\lesssim T^{\frac14} \frac{\epsilon}{\sqrt{\mu}} \left(  \| \| \zeta \|_{H^{5}_{x}} \| \bV \|_{H^{5}_{x}} \|_{L^{\frac43}_{t}(0,T)} +\|  \| \bV \|_{H^{5}_{x}}^2 \|_{L^{\frac43}_{t}(0,T)} \right)\\
&\lesssim T \frac{\epsilon}{\sqrt{\mu}} M^2.
\end{aligned}
\]

The first bound follows by H\"older's inequality. \colr One can similarly obtain the second bound. Note that by differentiating one time we can also get a bound on $\| \nabla \tfrac{\nabla}{\Delta} \cdot \bV \|_{L^{2}_{t}(0,T;W^{1,\infty}_{x}(\R^{2}))}$ which will be useful in the following.\colb

\begin{comment}
\br\label{notL1}
Note that one can not get bounds on $B_{1}$ when $\tilde r=\infty$ since in order to control $\frac{\nabla}{|D|} \cdot ((\bV \cdot \nabla) \bV)$ in $L^{1}_{x}$, $(\bV \cdot \nabla) \bV$ must have spatial integral zero (see for instance \cite[Corollary 2.4.8]{Grafakov_modernFourier} for a proof of this fact) which is not true in general.
\er
\end{comment}

On the other hand, since $\tilde \bV = \nabla^{\perp}\frac{\nabla^{\perp}}{\Delta} \cdot \tilde \bV$,
\[
\partial_{t} (\nabla^{\perp}\tfrac{\nabla^{\perp}}{\Delta} \cdot \bV - \tilde \bV) + \nabla^{\perp}\frac{\nabla^{\perp}}{\Delta} \cdot ( (\tilde \bV \cdot \nabla ) (\bV-\tilde \bV) ) + \nabla^{\perp}\frac{\nabla^{\perp}}{\Delta} \cdot ( (\bV-\tilde \bV) \cdot \nabla ) \bV )  = 0
\]
so that integrating by parts
\begin{align*}
\frac12 \frac{d}{dt} \left( \| \nabla^{\perp}\tfrac{\nabla^{\perp}}{\Delta} \cdot \bV - \tilde \bV \|^2_{L^{2}_{x}}\right) = &- \underbrace{\int_{\R^{2}} ( (\tilde \bV \cdot \nabla ) (\bV-\tilde \bV) ) \cdot (\nabla^{\perp}\tfrac{\nabla^{\perp}}{\Delta} \cdot \bV - \tilde \bV)}_{:= J_{1}}\\
&- \underbrace{\int_{\R^{2}} (( (\bV-\tilde \bV) \cdot \nabla ) \bV) \cdot (\nabla^{\perp}\tfrac{\nabla^{\perp}}{\Delta} \cdot \bV - \tilde \bV)}_{:= J_{2}}.
\end{align*}
Then
\[
J_{1} = \int_{\R^{2}} ( (\tilde \bV \cdot \nabla ) (\nabla^{\perp}\tfrac{\nabla^{\perp}}{\Delta} \cdot \bV-\tilde \bV) ) \cdot (\nabla^{\perp}\tfrac{\nabla^{\perp}}{\Delta} \cdot \bV - \tilde \bV) + \int_{\R^{2}} ( (\tilde \bV \cdot \nabla ) \nabla \tfrac{\nabla}{\Delta} \cdot \bV)  \cdot (\nabla^{\perp}\tfrac{\nabla^{\perp}}{\Delta} \cdot \bV - \tilde \bV)
\]
so that integrating by parts in the first integral and using that $\nabla \cdot \tilde \bV = 0$
\[
|J_{1}| \leq \| \nabla  \nabla \tfrac{\nabla}{\Delta} \cdot \bV \|_{L^{\infty}_{x}} \| \tilde \bV \|_{L^{2}_{x}} \| \nabla^{\perp}\tfrac{\nabla^{\perp}}{\Delta} \cdot \bV - \tilde \bV \|_{L^{2}_{x}}
\]
whereas
\[
|J_{2}| \leq \| \nabla \bV \|_{L^{\infty}_{x}} \| \nabla^{\perp}\tfrac{\nabla^{\perp}}{\Delta} \cdot \bV - \tilde \bV \|^2_{L^{2}_{x}} + \| \nabla \tfrac{\nabla}{\Delta} \cdot \bV \|_{L^{\infty}_{x}} \| \nabla \bV \|_{L^{2}_{x}} \| \nabla^{\perp}\tfrac{\nabla^{\perp}}{\Delta} \cdot \bV - \tilde \bV \|_{L^{2}_{x}}
\]
which yields
\[
\colp \frac{d}{dt} \left( \| \nabla^{\perp}\tfrac{\nabla^{\perp}}{\Delta} \cdot \bV - \tilde \bV \|^2_{L^{2}_{x}}\right) \leq C  M  \| \nabla^{\perp}\tfrac{\nabla^{\perp}}{\Delta} \cdot \bV - \tilde \bV \|^2_{L^{2}_{x}} + M  \| \nabla \tfrac{\nabla}{\Delta} \cdot \bV \|_{W^{1,\infty}_{x}} \| \nabla^{\perp}\tfrac{\nabla^{\perp}}{\Delta} \cdot \bV - \tilde \bV \|_{L^{2}_{x}}.\colb
\]
Bounds on $\nabla^{\perp}\tfrac{\nabla^{\perp}}{\Delta} \cdot \bV - \tilde \bV$ follow from Gr\"onwall's inequality and previous bounds obtained on $\| \nabla \tfrac{\nabla}{\Delta} \cdot \bV \|_{L^{2}_{t}(0,T;W^{1,\infty}_{x})}$.

Concerning Theorem \ref{rigidlid2d2}, one can use the previous strategy together with \eqref{strichartz1/22d}\footnote{\colr It is not necessary to split the low and the high frequencies in that case since we do not use Proposition \ref{propstrichartz12d}.\colb}. The last bound is a consequence of Morawetz-type estimates established in Proposition \ref{morawetz}.

\section{Other abcd Boussinesq systems}\label{s:Boussiabcd}

In the previous two sections we chose to present the rigid lid limit on one specific Boussinesq system (the case $a=b=c=0$ and $d=\tfrac13$). There are other abcd-Boussinesq systems
\be\label{abcd_boussinesq_eq}
\left\{
\begin{aligned}
&\epsilon \left(1- \mu b \Delta \right) \partial_{t} \zeta + \nabla\cdot \left( \left[1 + \epsilon \zeta \right] \bV \right) + \mu a \Delta \nabla \cdot \bV = 0,\\
&\epsilon \left(1- \mu d \nabla \nabla \cdot \right) \partial_{t} \bV + \nabla \zeta + \epsilon (\bV \cdot \nabla) \bV + \mu c \Delta \nabla \zeta = 0.
\end{aligned}
\right.
\ee

In the following, we assume that 
\be\label{condition_abcd}
b \geq 0 \text{  ,  } d \geq 0 \text{  ,  } a \leq 0 \text{  ,  } c \leq 0
\ee
in order to get the wellposedness of the system (see for instance \cite{bona_chen_saut_derivation}). 

We introduce 
\[
g(y) = y \sqrt{\frac{(1- a y^2)(1- c y^2)}{(1+ b y^2)(1+ d y^2)}} \text{   and   } R(y) = \sqrt{\frac{(1- a y^2) (1+ d y^2)}{(1+ b y^{2}) (1- c y^{2})}}.
\]

As before, in the 1d case if we denote by $\bU = (\zeta, V)^{T}$ we have the following system 
\bes
\epsilon \partial_{t} \bU + A(\partial_{x}) \bU = \epsilon  F(\zeta,V)
\ees
where
\begin{align*}
&A(\partial_{x}) = \begin{pmatrix} 0 & (1-\mu b \partial_{x}^{2})^{-1}(1+\mu a \partial_{x}^{2}) \partial_{x} \\ (1-\mu d \partial_{x}^{2})^{-1}(1+\mu c \partial_{x}^{2}) \partial_{x} & 0 \end{pmatrix}\\
&F(\zeta,V) = - \begin{pmatrix} (1 - \mu b \partial_{x}^{2})^{-1} \partial_{x} (\zeta V) \\ (1 - \mu d \partial_{x}^{2})^{-1} \partial_{x}( \frac{1}{2} V^{2} )\end{pmatrix}.
\end{align*}
Note that
\bes
\exp(tA(\partial_{x})) =  \begin{pmatrix} \cos (\frac{t}{\sqrt{\mu}} g(\sqrt{\mu} D)) & R(\sqrt{\mu} D) \sin (\frac{t}{\sqrt{\mu}} g(\sqrt{\mu} D)) \\ \frac{1}{R(\sqrt{\mu} D)} \sin (\frac{t}{\sqrt{\mu}} g(\sqrt{\mu} D)) & \cos (\frac{t}{\sqrt{\mu}} g(\sqrt{\mu} D)) \end{pmatrix}.
\ees

In the 2d case if we denote by $\bU = (\zeta,\nabla \cdot \bV)^{T}$, we get the following system
\bes
\epsilon \partial_{t} \bU + A(D) \bU = \epsilon  F(\zeta,\bV)
\ees
where
\begin{align*}
&A(D) = \begin{pmatrix} 0 & (1-\mu b \Delta)^{-1}(1+\mu a \Delta) \\ (1-\mu d \Delta)^{-1}(1+\mu c \Delta) \Delta & 0 \end{pmatrix}\\
& F(\zeta,\bV) = - \begin{pmatrix} (1 - \mu b \Delta)^{-1} \nabla \cdot (\zeta \bV) \\ (1 - \mu d \Delta)^{-1} \nabla \cdot ( (\bV \cdot \nabla ) \bV) \end{pmatrix}.
\end{align*}
Note that
\bes
\exp(tA(D)) = \begin{pmatrix} \cos (\frac{t}{\sqrt{\mu}} g(\sqrt{\mu} |D|)) & \frac{R(\sqrt{\mu} |D|)}{|D|} \sin (\frac{t}{\sqrt{\mu}} g(\sqrt{\mu} |D|)) \\ \frac{|D|}{R(|D|)} \sin (\frac{t}{\sqrt{\mu}} g(\sqrt{\mu} |D|)) & \cos (\frac{t}{\sqrt{\mu}} g(\sqrt{\mu} |D|)) \end{pmatrix}.
\ees

The strategy presented in the previous two sections together with ad hoc dispersive estimates provide similar results for System \eqref{abcd_boussinesq_eq}, with a \colr rate of convergence \colb depending on how dispersive System \eqref{abcd_boussinesq_eq} is. Existence of solutions of \eqref{abcd_boussinesq_eq} on an existence time independent of $\epsilon$ uniformly with respect to $\mu \in (0,1]$ can easily be adapted from \cite{saut_xu_system_boussi,saut_xu_system_boussi_II}.

The phase $g$ satisfies the following properties that are carefully studied in \cite[Section 3.5]{my_weaklydispestim}. Firstly if $a+b+c+d \neq 0$
\[
g'(r) -1 \underset{r \sim 0}{\sim}  - \frac{3(a+b+c+d)}{2} r^2 \text{   ,   } g''(r) \underset{r \sim 0}{\sim} - 3(a+b+c+d) r
\]
whereas if $a+b+c+d = 0$
\[
g'(r) - 1 \underset{r \sim 0}{\sim} - \frac{5(a+b)(b+c)}{2} r^4 \text{   ,   } g''(r) \underset{r \sim 0}{\sim} - 10(a+b)(b+c) r^3.
\]
Secondly there exists $\alpha \in [-6,1] \cap \Z$, $\ell, \Lambda_{1},\Lambda_{2} \in \R$ such that
\[
g'(r) - \ell \underset{\infty}{\sim} \Gamma_{1} r^{\alpha+1} \text{  ,  } g''(r) \underset{\infty}{\sim} (\alpha+1) \Gamma_{1} r^{\alpha} \text{  ,  } g'''(r) \leq \Gamma_{2} r^{\alpha-1}.
\]
The exact value of $\alpha$ and $\ell$ (that depends on $a,b,c,d$) can be found in \cite[Table 1]{my_weaklydispestim}.
Finally one can prove that $|g'|+|g''|+|g'''|>0$ on $\R^{+}$ (see \cite[Lemma 3.4]{my_weaklydispestim}).

We can now state our results. We begin with the case $n=1$ in the case $a+b+c+d \neq 0$.

\bt\label{rigidlid1dabcdneq0}
Let $n=1$. Let $a,b,c,d$ satisfying \eqref{condition_abcd} and $a+b+c+d \neq 0$. Let $M>0$, $T>0$, $\epsilon \in (0,1]$ and $\mu \in (0,1]$. Let $(\zeta,V) \in \mathcal{C}([0,T];(H^{3} \times H^{3})(\R))$ a solution of \eqref{abcd_boussinesq_eq} with initial datum $(\zeta_{0},V_{0})$ such that
\[
\| (\zeta,V) \|_{L^{\infty}(0,T;H^{3}(\R) \times H^{3}(\R))} \leq M.
\]
There exists $p \in \N$ with $p \geq 3$ and a constant \colr $C>0$ depending only on $p$ and $a,b,c,d$ \colb such that for any $q,r \geq 2$ with $\frac{1}{q} + \frac{1}{p r} = \frac{1}{2p}$
\[
\begin{aligned}
&\left\| \bp \zeta \\ V \ep - e^{-\frac{t}{\epsilon} A(\partial_{x})} \bp \zeta_{0} \\ V_{0} \ep \right\|_{L_{t}^{q}(0,T;L_{x}^{r}(\R))} \leq \left( \frac{\epsilon}{\mu} \right)^{\frac{1}{2 p} + \frac{1}{q}} \colr M^2 \colb T^{\frac{2 p -1}{2 p}} C,\\
&\left\| \bp \zeta \\ V \ep \right\|_{L_{t}^{q}(0,T;L_{x}^{r}(\R))} \leq \left( \frac{\epsilon}{\mu} \right)^{\frac{1}{q}} \left( \colr M \colb + \colr M^2 \colb T^{\frac{2 p -1}{2p}} \left( \frac{\epsilon}{\mu} \right)^{\frac{1}{2 p}} \right) C.
\end{aligned}
\]
Furthermore if  $|g''|>0$ on $\R^{\ast}_{+}$, for any $\tilde q, \tilde r \geq 2$ with $\frac{1}{\tilde q} + \frac{1}{2 \tilde r} = \frac14$ 
\[
\left\| \bp \zeta \\ V \ep - e^{-\frac{t}{\epsilon} A(\partial_{x})} \bp \zeta_{0} \\ V_{0} \ep \right\|_{L_{t}^{\tilde q}(0,T;L_{x}^{\tilde r}(\R))} \leq \left( \frac{\epsilon}{\mu} \right)^{\frac14 + \frac{1}{\tilde q}} \colr M^2 \colb T^{\frac34} C.
\]
Finally denoting $p_{0}=1$ if $|g'|>0$ on $\R^{+}$, $p_{0}=2$ if $|g'|+|g''|>0$ on $\R^{+}$ and $p_{0}=3$ otherwise, we have
\[
\sup_{x_{0} \in \R} \left\| e^{-(x-x_{0})^2} \bp \zeta \\ V \ep \right\|_{L_{t}^{2}(0,T;L^{2}_{x}(\R))} \leq \epsilon^{\frac{1}{2p_{0}}}( \colr M \colb + \colr M^2 \colb T) C.
\]
\et 

If we denote by $m$ the maximum among the multiplicities of positive zeros of $g''$ then one can take $p=\max(m+2,3)$. The proof of the previous theorem follows from dispersive estimates based on Lemma \ref{decay1d}\footnote{\colr with $\beta=1,s=0,l=p$ and $\beta=1,s=\frac12,l=2$ for the third estimate.\colb} and the properties of $g$. We only provide a proof of the last point. Let $\chi_{1}$ be a smooth bounded function supported on $\{ |g'| > 0 \}$ and $\chi_{2}$ a smooth compactly supported function supported on $\{ |g''|+ |g'''| > 0 \}$ with $0 \notin \text{supp}(\chi_{2})$. On one hand we get from Proposition \ref{morawetz} that
\[
\sup_{x_{0} \in \R} \left\| e^{-(x-x_{0})^2}  \chi_{1}(\sqrt{\mu} |D|) \bp \zeta \\ V \ep \right\|_{L_{t}^{2}(0,T;L^{2}_{x}(\R))} \leq \epsilon^{\frac12}(\colr M \colb + \colr M^2 \colb T) C.
\]
On the other hand using Lemma \ref{decay1d}(ii)\footnote{\colr with $l=p_{0}$, $\alpha=0$.\colb} together with Bernstein's Lemma \ref{Bernstein} we have
\[
\left\| e^{\pm \frac{\iD t}{\sqrt{\mu}} g(\sqrt{\mu} |D|)} \chi_{2}(\sqrt{\mu} |D|)) f \right\|_{L^{\infty}} \lesssim \frac{1}{|t|^{\frac{1}{p_{0}}}} \| |D|^{\frac{p_{0}-1}{p_{0}}} \chi_{2}(\sqrt{\mu} |D|)) f \|_{L^{1}}
\]
so that from corresponding Strichartz estimates we obtain
\[
\left\| \chi_{2}(\sqrt{\mu} |D|)) \bp \zeta \\ V \ep \right\|_{L_{t}^{2 p_{0}}(0,T;L_{x}^{\infty}(\R))} \leq\epsilon^{\frac{1}{2p_{0}}} \left( \colr M \colb + \colr M^2 \colb T^{\frac{2 p_{0} -1}{2 p_{0}}} \epsilon^{\frac{1}{2 p_{0}}} \right) C.
\]
The third point of the theorem follows from the fact that $|g'|+|g''|+|g'''|>0$ on $\R^{+}$ and that there exists $0<y_{0} \leq y_{1}$ such that $g'>0$ on $[0,y_{0}] \cup [y_{1},\infty)$.

\br As noted in Remark \ref{remarkKdV}  concerning Theorem \ref{rigidlid1d}, in the case $\epsilon \sim \mu$ as in \cite{bona_chen_saut_derivation,bona_colin_lannes} or when \colr $\mu = \mathcal{O}(\epsilon)$\colb, the first estimate of Theorem \ref{rigidlid1dabcdneq0} does not provide a convergence result as $\epsilon \to 0$ so that nonlinear terms must be taken into account and that asymptotic models like a system of decoupling KdV equations 
\[
\epsilon \partial_{t} g_{\pm} \pm \partial_{x} g_{\pm} \pm \mu \frac{a+b+c+d}{6} \partial_{x}^{3} g_{\pm} \pm \frac32 \epsilon g_{\pm} \partial_{x} g_{\pm} = 0
\]
becomes relevant. A proof of such a result can be adapted from for instance \cite[Section 7.3.2]{Lannes_ww} together with the symmetrizers and energy estimates from \cite{saut_xu_system_boussi,saut_xu_system_boussi_II}.
\er

We now consider the case $n=1$ in the case $a+b+c+d = 0$. We introduce the condition
\be\label{condition_abcd2}
((a+b)(a+d)(c+b)(c+d))^2+(a+b+c+d)^2 > 0
\ee
which avoids the situation where $g(r) \equiv r $ that provides a non dispersive system when $n=1$.

\bt\label{rigidlid1dabcd=0}
Let $n=1$. Let $a,b,c,d$ satisfying \eqref{condition_abcd} with $a+b+c+d=0$. Let $M>0$, $T>0$, $\epsilon \in (0,1]$ and $\mu \in (0,1]$. Let $(\zeta,V) \in \mathcal{C}([0,T];(H^{3} \times H^{3})(\R))$ a solution of \eqref{abcd_boussinesq_eq} with initial datum $(\zeta_{0},V_{0})$ such that
\[
\| (\zeta,V) \|_{L^{\infty}(0,T;H^{3}(\R) \times H^{3}(\R))} \leq M.
\]
If $a,b,c,d$ satisfy \eqref{condition_abcd2}, there exists $p \in \N$ with $p \geq 5$ and a constant \colr $C>0$ depending only on $p$ and $a,b,c,d$ \colb such that for any $q,r \geq 2$ with $\frac{1}{q} + \frac{1}{p r} = \frac{1}{2p}$
\[
\begin{aligned}
&\left\| \bp \zeta \\ V \ep - e^{-\frac{t}{\epsilon} A(\partial_{x})} \bp \zeta_{0} \\ V_{0} \ep \right\|_{L_{t}^{q}(0,T;L_{x}^{r}(\R))} \leq \left( \frac{\epsilon}{\mu^{2}} \right)^{\frac{1}{2 p} + \frac{1}{q}} \colr M^2 \colb T^{\frac{2 p -1}{2 p}} C,\\
&\left\| \bp \zeta \\ V \ep \right\|_{L_{t}^{q}(0,T;L_{x}^{r}(\R))} \leq \left( \frac{\epsilon}{\mu^{2}} \right)^{\frac{1}{q}} \left( \colr M \colb + \colr M^2 \colb T^{\frac{2 p -1}{2p}} \left( \frac{\epsilon}{\mu^{2}} \right)^{\frac{1}{2 p}} \right) C.
\end{aligned}
\]
Furthermore if some $l \in \{2,3,4\}$ we have $\sum_{k=2}^{l} |g^{(k)}|>0$ on $\R^{\ast}_{+}$ and if we denote by $\sigma = \min(\frac25,\frac{1}{l})$, for any $\tilde q, \tilde r \geq 2$ with $\frac{1}{\tilde q} +  \frac{\sigma}{\tilde r} = \frac{\sigma}{2}$ 
\[
\left\| \bp \zeta \\ V \ep - e^{-\frac{t}{\epsilon} A(\partial_{x})} \bp \zeta_{0} \\ V_{0} \ep \right\|_{L_{t}^{\tilde q}(0,T;L_{x}^{\tilde r}(\R))} \leq \left( \frac{\epsilon}{\mu^{2}} \right)^{\frac{\sigma}{2} + \frac{1}{\tilde q}} \colr M^2 \colb T^{\frac{2-\sigma}{2}} C.
\]
Finally denoting $p_{0}=1$ if $|g'|>0$ on $\R^{+}$, $p_{0}=2$ if $|g'|+|g''|>0$ on $\R^{+}$ and $p_{0}=3$ otherwise, we have
\[
\sup_{x_{0} \in \R} \left\| e^{-(x-x_{0})^2} \bp \zeta \\ V \ep \right\|_{L_{t}^{2}(0,T;L^{2}_{x}(\R))} \leq \epsilon^{\frac{1}{2p_{0}}}( \colr M \colb + \colr M^2 \colb T) C.
\]
\et 

If we denote by $m$ the maximum among the multiplicities of positive zeros of $g''$ then one can take $p=\max(m+2,5)$.  Again one can obtain dispersive estimates thanks to Lemma \ref{decay1d}\footnote{\colr with $\beta=3,s=0,l=p$ and $\beta=3,s=5 \sigma -1$ for the third estimate.\colb} and the previous properties on the phase $g$. Note that the ratio $\frac{\epsilon}{\mu^{2}}$ comes from low frequency estimates: if $\chi$ is a smooth compactly supported function whose support is small enough and that is equal to $1$ near $0$ and if $a+b+c+d=0$
\[
\begin{aligned}
&\left\| e^{\pm \frac{\iD t}{\epsilon \sqrt{\mu}} g(\sqrt{\mu} |D|)} \chi(\sqrt{\mu} |D|) \partial_{x} f \right\|_{L^{\infty}} \lesssim \frac{1}{|t|^{\frac25}} \left( \frac{\epsilon}{\mu^{2}} \right)^{\frac25} \| f \|_{L^{1}},\\
&\left\| e^{\pm \frac{\iD t}{\epsilon \sqrt{\mu}} g(\sqrt{\mu} |D|)} \chi(\sqrt{\mu} |D|)  f \right\|_{L^{\infty}} \lesssim \frac{1}{|t|^{\frac15}} \left( \frac{\epsilon}{\mu^{2}} \right)^{\frac15} \| f \|_{L^{1}}.
\end{aligned}
\]

\br When $\epsilon \sim \mu$ as in \cite{bona_chen_saut_derivation,bona_colin_lannes} or when \colr$\mu = \mathcal{O}(\epsilon)$\colb, the first estimate of Theorem \ref{rigidlid1dabcdneq0} does not provide a convergence result as $\epsilon \to 0$. Again nonlinear terms must be taken into account. Note however that here one must consider a system of decoupling Burgers equations 
\[
\epsilon \partial_{t} g_{\pm} \pm \partial_{x} g_{\pm} \pm \frac32 \epsilon g_{\pm} \partial_{x} g_{\pm} = 0.
\]
A proof of such a result can be adapted from for instance \cite[Section 7.3.2]{Lannes_ww} together with the symmetrizers and energy estimates from \cite{saut_xu_system_boussi,saut_xu_system_boussi_II}.
\er

We now consider the case $n=2$ with $a+b+c+d \neq 0$.

\bt\label{rigidlid2dabcdneq0}
Let $n=2$. Let $a,b,c,d$ satisfying \eqref{condition_abcd} with $a+b+c+d \neq 0$. Let $M>0$, $T>0$, $\epsilon \in (0,1]$ and $\mu \in (0,1]$.  Let $(\zeta,\bV) \in \mathcal{C}([0,T];(H^{6} \times H^{6})(\R^2))$ be a solution of \eqref{abcd_boussinesq_eq} with initial datum $(\zeta_{0},\bV_{0})$, let $\tilde \bV \in \mathcal{C}([0,T];L^{2}(\R^2))$ be a solution of the incompressible Euler equation  \eqref{incompEuler} with initial datum $\nabla^{\perp} \frac{\nabla^{\perp}}{\Delta} \cdot \bV_{0}$ such that
\[
\| (\zeta,\bV) \|_{L^{\infty}(0,T;(H^{6} \times H^{6})(\R^2))} + \| \tilde \bV \|_{L^{\infty}(0,T;L^{2}(\R^2))} \leq M.
\]

If $g'$ and $g''$ do not vanish on $\R^{\ast}_{+}$, $(\zeta,\bV)$ satisfy the same estimates as in Theorem \ref{rigidlid2d1}.

If $g'$ does not vanish on $\R^{+}$ but $g''$ vanishes on $\R^{\ast}_{+}$, there exists $\sigma \in (\tfrac12,1)$ and a constant \colr$C>0$ depending only on $a,b,c,d$ \colb such that for any $q,r \geq 2$ with $\frac{1}{q} + \frac{\sigma}{r} = \frac{\sigma}{2}$
\[
\begin{aligned}
&\left\| \bp \zeta \\ \frac{\nabla \nabla}{\Delta} \cdot \bV \ep - \bp 1 & 0 \\ 0 & \frac{\partial_{1}}{\Delta} \\ 0 & \frac{\partial_{2}}{\Delta}  \ep \exp(-\tfrac{t}{\epsilon}A(D)) \bp \zeta_{0} \\ \nabla \cdot \bV_{0} \ep \right\|_{L_{t}^{q}(0,T;L_{x}^{r}(\R^2))}  \leq \left( \frac{\epsilon}{\sqrt{\mu}} \right)^{\colr \frac{\sigma}{2} \colb + \frac{1}{q}} \colr M^2 \colb T^{\frac{2-\sigma}{2}} C,\\
&\left\|  \bp \zeta \\ \frac{\nabla \nabla }{\Delta} \cdot \bV \ep  \right\|_{L_{t}^{q}(0,T;L_{x}^{r}(\R^2))}  \leq \left( \frac{\epsilon}{\sqrt{\mu}} \right)^{\frac{1}{q}} \left( \colr M \colb + \colr M^2 \colb T^{\frac{2-\sigma}{2}} \left( \frac{\epsilon}{\sqrt{\mu}} \right)^{\colr \frac{\sigma}{2} \colb} \right) C,\\
&\colp\left\| \nabla^{\perp} \frac{\nabla^{\perp}}{\Delta}  \cdot \bV - \tilde \bV \right\|_{L_{t}^{\infty}(0,T;L_{x}^{2}(\R^2))} \leq \left( \frac{\epsilon}{\sqrt{\mu}} \right)^{ \frac{\sigma}{2} } \left(  M  +  M^2  T^{\frac{2-\sigma}{\sigma}} \left( \frac{\epsilon}{\sqrt{\mu}} \right)^{ \frac{\sigma}{2} } \right) M T^{\frac{2-\sigma}{2}}  e^{C  M  T} C.\colb
\end{aligned}
\]
Furthermore, let $p=2$ if $|g'|+|g''|>0$ on $\R^{+}$ and $p=3$ otherwise, there exists a constant \colr $C>0$ depending only on $a,b,c,d$ \colb such that for any $\tilde q, \tilde r \geq 2$ with $\frac{1}{\tilde q} + \frac{1}{p \tilde r} = \frac{1}{2p}$
\[
\begin{aligned}
&\left\| \bp \zeta \\ \frac{\nabla \nabla}{\Delta} \cdot \bV \ep - \bp 1 & 0 \\ 0 & \frac{\partial_{1}}{\Delta} \\ 0 & \frac{\partial_{2}}{\Delta}  \ep \exp(-\tfrac{t}{\epsilon}A(D)) \bp \zeta_{0} \\ \nabla \cdot \bV_{0} \ep \right\|_{L_{t}^{\tilde q}(0,T;L_{x}^{\tilde r}(\R^2))} \leq \epsilon^{\frac{1}{2 p} + \frac{1}{\tilde q}} \colr M^2 \colb T^{\frac{2 p-1}{2 p}} C,\\
&\left\| \bp \zeta \\ \frac{\nabla \nabla }{\Delta} \cdot \bV \ep \right\|_{L_{t}^{\tilde q}(0,T;L_{x}^{\tilde r}(\R^2))} \leq \epsilon^{\frac{1}{\tilde q}} \left( \colr M \colb + \colr M^2 \colb T^{\frac{2 p-1}{2 p}} \epsilon^{\frac{1}{2 p}}  \right) C,\\
&\colp \left\| \nabla^{\perp} \frac{\nabla^{\perp}}{\Delta}  \cdot \bV - \tilde \bV \right\|_{L_{t}^{\infty}(0,T;L_{x}^{2}(\R^2))} \leq \epsilon^{\frac{1}{2p}} \left( M + M^2 T^{\frac{2p-1}{2p}} \epsilon^{\frac{1}{2p}} \right) M T^{\frac{2p-1}{2p}} e^{C M T} C. \colb
\end{aligned}
\]
Finally denoting $p_{0}=1$ if $|g'|>0$ on $\R^{+}$, $p_{0}=2$ if $|g'|+|g''|>0$ on $\R^{+}$ and $p_{0}=3$ otherwise, we have
\[
\sup_{x_{0} \in \R^2} \left\| e^{-(x-x_{0})^2} \bp \zeta \\ \nabla \frac{\nabla}{\Delta} \cdot \bV \ep \right\|_{L_{t}^{2}(0,T;L^{2}_{x}(\R^2))} \leq \epsilon^{\frac{1}{2p_{0}}}(\colr M \colb + \colr M^2 \colb T) C.
\]
\et 
If we denote by $m$ the maximum among the multiplicities of positive zeros of $g''$ then one can take $\sigma=\frac{m+4}{2m+4}$. Again the key point are dispersive estimates that can be obtained from Lemmas \ref{decay2d}\footnote{\colr with $\beta=1$, $s=\alpha$ if $\ell=0$ and $\beta=1$, $s=\frac{\alpha-1}{2}$ if $\ell \neq 0$. \colb}, \ref{decay2d1/2} and \ref{morawetz}. 

Finally a similar result can be obtained in the case $n=2$ with $a+b+c+d = 0$. 

\bt\label{rigidlid2dabcd=0}
Let $n=2$. Let $a,b,c,d$ satisfying \eqref{condition_abcd} with $a+b+c+d = 0$. Let $M>0$, $T>0$, $\epsilon \in (0,1]$ and $\mu \in (0,1]$.  Let $(\zeta,\bV) \in \mathcal{C}([0,T];(H^{6} \times H^{6})(\R^2))$ be a solution of \eqref{abcd_boussinesq_eq} with initial datum $(\zeta_{0},\bV_{0})$, let $\tilde \bV \in \mathcal{C}([0,T];L^{2}(\R^2))$ be a solution of the incompressible Euler equation  \eqref{incompEuler} with initial datum $\nabla^{\perp} \frac{\nabla^{\perp}}{\Delta} \cdot \bV_{0}$ such that
\[
\| (\zeta,\bV) \|_{L^{\infty}(0,T;(H^{6} \times H^{6})(\R^2))} + \| \tilde \bV \|_{L^{\infty}(0,T;L^{2}(\R^2))} \leq M.
\]

If $a,b,c,d$ satisfies \eqref{condition_abcd2} and $g'>0$ on $\R^{+}$, there exists $\sigma \in (\tfrac12,\tfrac45]$ and \colr$C>0$ depending only on $a,b,c,d$\colb, such that for any $q,r \geq 2$ with $\frac{1}{q} + \frac{\sigma}{r} = \frac{\sigma}{2}$
\[
\begin{aligned}
&\left\|  \bp \zeta \\ \frac{\nabla \nabla}{\Delta} \cdot \bV \ep - \bp 1 & 0 \\ 0 & \frac{\partial_{1}}{\Delta} \\ 0 & \frac{\partial_{2}}{\Delta}  \ep \exp(-\tfrac{t}{\epsilon}A(D)) \bp \zeta_{0} \\ \nabla \cdot \bV_{0} \ep \right\|_{L_{t}^{q}(0,T;L_{x}^{r}(\R^2))}  \leq \left( \frac{\epsilon}{\mu^{\frac34}} \right)^{\colr \frac{\sigma}{2} \colb + \frac{1}{q}} \colr M^2 \colb T^{\frac{2-\sigma}{2}} C,\\
&\left\| \bp \zeta \\ \frac{\nabla \nabla }{\Delta} \cdot \bV \ep \right\|_{L_{t}^{q}(0,T;L_{x}^{r}(\R^2))}  \leq \left( \frac{\epsilon}{\mu^{\frac34}} \right)^{\frac{1}{q}} \left( \colr M \colb + \colr M^2 \colb T^{\frac{2-\sigma}{2}} \left( \frac{\epsilon}{\mu^{\frac34}} \right)^{\colr \frac{\sigma}{2} \colb} \right) C,\\
&\colp \left\| \nabla^{\perp} \frac{\nabla^{\perp}}{\Delta}  \cdot \bV - \tilde \bV \right\|_{L_{t}^{\infty}(0,T;L_{x}^{2}(\R^2))} \leq \left( \frac{\epsilon}{\mu^{\frac34}} \right)^{ \frac{\sigma}{2} } \left(  M +   M^2  T^{\frac{2-\sigma}{2}} \left( \frac{\epsilon}{\mu^{\frac34}} \right)^{ \frac{\sigma}{2} } \right) M T^{\frac{2-\sigma}{2}}  e^{C  M  T} C.\colb
\end{aligned}
\]
Furthermore, let $p=2$ if $|g'|+|g''|>0$ on $\R^{+}$ and $p=3$ otherwise, there exists a constant \colr $C>0$ depending only on $a,b,c,d$ \colb such that for any $\tilde q, \tilde r \geq 2$ with $\frac{1}{\tilde q} + \frac{1}{p \tilde r} = \frac{1}{2p}$
\[
\begin{aligned}
&\left\| \bp \zeta \\ \frac{\nabla \nabla}{\Delta} \cdot \bV \ep - \bp 1 & 0 \\ 0 & \frac{\partial_{1}}{\Delta} \\ 0 & \frac{\partial_{2}}{\Delta}  \ep \exp(-\tfrac{t}{\epsilon}A(D)) \bp \zeta_{0} \\ \nabla \cdot \bV_{0} \ep \right\|_{L_{t}^{\tilde q}(0,T;L_{x}^{\tilde r}(\R^2))} \leq \epsilon^{\frac{1}{2p} + \frac{1}{\tilde q}} \colr M^2 \colb T^{\frac{2p-1}{2p}} C,\\
&\left\| \bp \zeta \\ \frac{\nabla \nabla }{\Delta} \cdot \bV \ep \right\|_{L_{t}^{\tilde q}(0,T;L_{x}^{\tilde r}(\R^2))} \leq \epsilon^{\frac{1}{\tilde q}} \left( \colr M \colb +  \colr M^2 \colb T^{\frac{2p-1}{2p}} \epsilon^{\frac{1}{2p}} \right) C,\\
&\colp \left\| \nabla^{\perp} \frac{\nabla^{\perp}}{\Delta}  \cdot \bV - \tilde \bV \right\|_{L_{t}^{\infty}(0,T;L_{x}^{2}(\R^2))} \leq \epsilon^{\frac{1}{2p}} \left( M + M^2 T^{\frac{2p-1}{2p}} \epsilon^{\frac{1}{2p}} \right) M T^{\frac{2p-1}{2p}} e^{C M T} C.\colb
\end{aligned}
\]
Finally denoting $p_{0}=1$ if $|g'|>0$ on $\R^{+}$, $p_{0}=2$ if $|g'|+|g''|>0$ on $\R^{+}$ and $p_{0}=3$ otherwise, we have
\[
\sup_{x_{0} \in \R^2} \left\| e^{-(x-x_{0})^2} \bp \zeta \\ \nabla \frac{\nabla}{\Delta} \cdot \bV \ep \right\|_{L_{t}^{2}(0,T;L^{2}_{x}(\R^2))} \leq \epsilon^{\frac{1}{2p_{0}}}(\colr M \colb + \colr M^2 \colb T) C.
\]
\et 
If we denote by $m$ the maximum among the multiplicities of positive zeros of $g''$ then one can take $\sigma=\min(\frac{m+4}{2m+4},\tfrac45)$.  Again the key point are dispersive estimates that can be obtained from Lemmas \ref{decay2d}\footnote{\colr with $\beta=3$, $s=\alpha$ if $\ell=0$ and $\beta=3$, $s=\frac{\alpha-1}{2}$ if $\ell \neq 0$. \colb}, \ref{decay2d1/2} and \ref{morawetz}.

\section{The Green-Naghdi equations}\label{s:GN}

The Green-Naghdi equations read as

\be\label{GN_eq}
\left\{
\begin{aligned}
&\epsilon \partial_{t} \zeta + \nabla\cdot \left( \left[1 + \epsilon \zeta \right] \bV \right) = 0,\\
&\epsilon (1+\mu \mathcal{T}[\epsilon \zeta]) \partial_{t} \bV + \nabla \zeta + \epsilon (\bV \cdot \nabla) \bV + \epsilon \mu \mathcal{Q}[\epsilon \zeta](\bV) = 0
\end{aligned}
\right.
\ee
where
\[
\begin{aligned}
&\mathcal{T}[\epsilon \zeta] \bW = - \frac{1}{3(1+\epsilon \zeta)} \nabla \left[ (1+\epsilon \zeta)^{3} \nabla \cdot \bW \right] \\
&\mathcal{Q}[\epsilon \zeta](\bV) =  - \frac{1}{3(1+\epsilon \zeta)} \nabla \left[ (1+\epsilon \zeta)^{3} ( (\bV \cdot \nabla) (\nabla \cdot \bV) - (\nabla \cdot \bV)^{2} ) \right].
\end{aligned}
\]
As before, in the 1d case if we denote by $\bU = (\zeta, V)^{T}$ we have the following system 
\bes
\epsilon \partial_{t} \bU + A(\partial_{x}) \bU = \epsilon  F(\zeta,V)
\ees
where
\begin{align*}
&A(\partial_{x}) = \begin{pmatrix} 0 & \partial_{x} \\ (1- \frac{\mu}{3} \partial_{x}^{2})^{-1} \partial_{x} & 0 \end{pmatrix}\\
&F(\zeta,V) = - \bp  \partial_{x} (\zeta V) \\ (1- \frac{\mu}{3} \partial_{x}^{2})^{-1} \left( \partial_{x} ( \tfrac12 V^{2} ) + \mu \mathcal{Q}[\epsilon \zeta](V) + \mu \mathcal{T}[\epsilon \zeta] \partial_{t} V + \frac{\mu}{3} \partial_{x}^{2} \partial_{t} V \right) \ep ,
\end{align*}
whereas in the 2d case if we denote by $\bU = (\zeta,\nabla \cdot \bV)^{T}$, we get 
\bes
\epsilon \partial_{t} \bU + A(D) \bU = \epsilon  F(\zeta,\bV)
\ees
where
\begin{align*}
&A(D) = \begin{pmatrix} 0 & 1 \\ (1- \frac{\mu}{3} \Delta)^{-1} \Delta & 0 \end{pmatrix} \\
& F(\zeta,\bV) = - \bp  \nabla \cdot (\zeta \bV) \\ (1- \frac{\mu}{3} \Delta)^{-1} \nabla \cdot \left( ( \bV \cdot \nabla ) \bV  + \mu \mathcal{Q}[\epsilon \zeta](\bV) + \mu \mathcal{T}[\epsilon \zeta] \partial_{t} \bV + \frac{\mu}{3} \Delta \partial_{t} \bV \right) \ep.
\end{align*}

We refer to \cite{israwi_green_naghdi,duchene_israwi_green_naghdi} (see also \cite{Li_green_naghdi,Fujiwara_Iguchi_2015}) for the existence of solutions of \eqref{GN_eq} on an existence time independent of $\epsilon$ that is uniform with respect to $\mu \in (0,1]$.

We consider now the case $n=1$.

\bt\label{rigidlid1dGN}
Let $M>0$, $T>0$, $h_{0} >0$, $\epsilon \in (0,1]$ and $\mu \in (0,1]$. Let $(\zeta,V) \in \mathcal{C}([0,T];(H^{5} \times H^{5})(\R))$ a solution of \eqref{GN_eq} with initial datum $(\zeta_{0},V_{0})$ such that
\[
\| (\zeta,V) \|_{L^{\infty}(0,T;H^{5}(\R) \times H^{5}(\R))} \leq M \text{   and   } 1+ \epsilon \zeta \geq h_{0} \text{   on   } [0,T].
\]
There exists a constant \colr$C_{d}>0$ polynomial in $M$ and $1/h_{0}$ \colb such that for any $q,r \geq 2$ with $\frac{1}{q} + \frac{1}{2r} = \frac{1}{6}$
\[
\begin{aligned}
&\left\| \bp \zeta \\ V \ep - e^{-\frac{t}{\epsilon} A(\partial_{x})} \bp \zeta_{0} \\ V_{0} \ep \right\|_{L_{t}^{q}(0,T;L_{x}^{r}(\R))} \leq \left( \frac{\epsilon}{\mu} \right)^{\frac16 + \frac{1}{q}} T^{\frac56} \colr C_{d}\colb,\\
&\left\| \bp \zeta \\ V \ep \right\|_{L_{t}^{q}(0,T;L_{x}^{r}(\R))} \leq \left( \frac{\epsilon}{\mu} \right)^{\frac{1}{q}} \left( 1+T^{\frac56} \left( \frac{\epsilon}{\mu} \right)^{\frac{1}{6}} \right) \colr C_{d}\colb,\\
&\sup_{x_{0} \in \R} \left\| e^{-(x-x_{0})^2} \bp \zeta \\ V \ep \right\|_{L_{t}^{2}(0,T;L^{2}_{x}(\R))} \leq \epsilon^{\frac12}(1+T) \colr C_{d}\colb.
\end{aligned}
\]
\et 
The proof follows from the same strategy as the proof of Theorem \ref{rigidlid1d} \colr together with \eqref{strichartz1/31d} and \eqref{retardedstrichartz1/31d}\footnote{\colr Note that the source term $F$ is not a derivative here so that one can not use \eqref{retardedstrichartz1/21d}.\colb}\colb. One must control two new terms. Standard product estimates provide
\begin{align*}
&\left\| \mathcal{Q}[\epsilon \zeta](\bV) \right\|_{W^{2,1}_{x}} \lesssim C \left( \tfrac{1}{h_{0}} , \| \zeta \|_{H^{3}} \right) \| V \|_{H^{5}}^{2}\\
&\left\| \mu \mathcal{T}[\epsilon \zeta] \partial_{t} V + \frac{\mu}{3} \partial_{x}^{2} \partial_{t} V \right\|_{W^{2,1}_{x}} \lesssim C \left( \tfrac{1}{h_{0}} , \| \zeta \|_{H^{3}} \right) \| \epsilon \mu \partial_{x} \partial_{t} V \|_{H^{3}} 
\end{align*}
and using for instance ideas from the proofs of \cite[Lemmas 1 and 2]{israwi_green_naghdi} and standard product estimates we obtain
\[
\begin{aligned}
\| \epsilon \mu \partial_{x} \partial_{t} V \|_{H^{3}} &\lesssim \epsilon C \left( \tfrac{1}{h_{0}} , \| \zeta \|_{H^{3}} \right) \left\| \nabla \zeta + \epsilon (V \cdot \nabla) V + \epsilon \mu \mathcal{Q}[\epsilon \zeta](V) \right\|_{H^{2}}\\
&\lesssim C \left( \tfrac{1}{h_{0}} , \| \zeta \|_{H^{3}} , \| V \|_{H^{5}} \right).
\end{aligned}
\]

\br As noted in Remark \ref{remarkKdV} concerning Theorem \ref{rigidlid1d}, in the case \colr$\mu = \mathcal{O}(\epsilon)$ \colb the first estimate of Theorem \ref{rigidlid1dGN} does not provide a convergence result as $\epsilon \to 0$ so that nonlinear terms must be taken into account and asymptotic models like a system of decoupling KdV equations or decoupling BBM equations become relevant. We refer to \cite[Chapter 7]{Lannes_ww}.
\er

We now consider the case $n=2$. Applying the operator $\nabla^{\perp} \cdot$ to the second equation of \eqref{GN_eq} and denoting by $\omega := \nabla^{\perp} \cdot \bV$ we get the following equation
\[
\partial_{t} \omega + (\bV \cdot \nabla) \omega + (\nabla \cdot \bV) \omega + \frac{\epsilon \mu \nabla \zeta^{\perp}}{3(1+\epsilon \zeta)^2} \cdot \nabla \left[ (1+\epsilon \zeta)^{3} \left( \nabla \cdot \partial_{t} \bV + (\bV \cdot \nabla) (\nabla \cdot \bV) - (\nabla \cdot \bV)^{2}  \right) \right] = 0.
\]

\bt\label{rigidlid2dGN}
Let $M>0$, $T>0$, $h_{0}>0$, $\epsilon \in (0,1]$ and $\mu \in (0,1]$. Let $(\zeta,\bV) \in \mathcal{C}([0,T];(H^{9} \times H^{9})(\R^2))$ be a solution of \eqref{GN_eq} with initial datum $(\zeta_{0},\bV_{0})$ and $\tilde \bV \in \mathcal{C}([0,T];L^{2}(\R^2))$ be a solution of the incompressible Euler equation \eqref{incompEuler} with initial datum $\nabla^{\perp} \frac{\nabla^{\perp}}{\Delta} \cdot \bV_{0}$ such that
\[
\| (\zeta,\bV) \|_{L^{\infty}(0,T;(H^{9} \times H^{9})(\R^2))} + \| \tilde \bV \|_{L^{\infty}(0,T;L^{2}(\R^2))} \leq M \text{   and   } 1+ \epsilon \zeta \geq h_{0} \text{   on   } [0,T].
\]
There exists a constant \colr$C_{d}>0$ polynomial in $M$ and $1/h_{0}$ \colb and a universal constant $C>0$ \colb such that for any $q,r \geq 2$ with $\frac{1}{q} + \frac{1}{r} = \frac{1}{2}$
\[
\begin{aligned}
&\left\| \bp \zeta \\ \frac{\nabla \nabla}{\Delta} \cdot \bV \ep - \bp 1 & 0 \\ 0 & \frac{\partial_{1}}{\Delta} \\ 0 & \frac{\partial_{2}}{\Delta}  \ep \exp(-\tfrac{t}{\epsilon}A(D)) \bp \zeta_{0} \\ \nabla \cdot \bV_{0} \ep \right\|_{L_{t}^{q}(0,T;L_{x}^{r}(\R^2))} \leq \left( \frac{\epsilon}{\sqrt{\mu}}  \ln(\colr 1 \colb + \tfrac{\mu}{\epsilon^2} T) \right)^{\frac12 + \frac{1}{q}} T^{\frac12} \colr C_{d}\colb\\
&\hspace{11cm} + \tfrac{\epsilon}{\sqrt{\mu}} T \colr C_{d}\colb,\\
&\left\| \bp \zeta \\ \frac{\nabla \nabla }{\Delta} \cdot \bV \ep   \right\|_{L_{t}^{q}(0,T;L_{x}^{r}(\R^2))} \leq \left( \frac{\epsilon}{\sqrt{\mu}}  \ln(\colr 1 \colb + \tfrac{\mu}{\epsilon^2} T) \right)^{\frac{1}{q}} \left( 1 +\left( \frac{\epsilon}{\sqrt{\mu}}  \ln(\colr 1 \colb + \tfrac{\mu}{\epsilon^2} T) T \right)^{\frac12} \right) \colr C_{d}\colb\\
&\hspace{5cm} + \tfrac{\epsilon}{\sqrt{\mu}} (1+T) \colr C_{d}\colb,\\
&\colp \left\| \nabla^{\perp} \frac{\nabla^{\perp}}{\Delta}  \cdot \bV - \tilde \bV \right\|_{L_{t}^{\infty}(0,T;L_{x}^{2}(\R^2))} \hspace{-0.5cm}  \leq \left( \frac{\epsilon}{\sqrt{\mu}} \ln( 1 + \tfrac{\mu}{\epsilon^2} T) \right)^{\frac12} \left( 1+  \left( \frac{\epsilon}{\sqrt{\mu}} \ln( 1  + \tfrac{\mu}{\epsilon^2} T) T \right)^{\frac12} \right) \sqrt{T} e^{ C M T} C_{d} \colb\\
&\hspace{6cm} \colp + \tfrac{\epsilon}{\sqrt{\mu}} (1+T) \sqrt{T} e^{ C M T} C_{d}\colb
\end{aligned}
\]
and there exists a constant \colr$\tilde C_{d}>0$ polynomial in $M$ and $1/h_{0}$ \colb and a universal constant $\tilde C>0$ \colb such that for any $\tilde q, \tilde r \geq 2$ with $\frac{1}{\tilde q} + \frac{1}{2 \tilde r} = \frac{1}{4}$
\[
\begin{aligned}
&\left\| \bp \zeta \\ \frac{\nabla \nabla}{\Delta} \cdot \bV \ep - \bp 1 & 0 \\ 0 & \frac{\partial_{1}}{\Delta} \\ 0 & \frac{\partial_{2}}{\Delta}  \ep \exp(-\tfrac{t}{\epsilon}A(D)) \bp \zeta_{0} \\ \nabla \cdot \bV_{0} \ep \right\|_{L_{t}^{\tilde q}(0,T;L_{x}^{\tilde r}(\R^2))} \leq \epsilon^{\frac14 + \frac{1}{\tilde q}} T^{\frac34} \colr \tilde C_{d}\colb,\\
&\left\| \bp \zeta \\ \frac{\nabla \nabla }{\Delta} \cdot \bV \ep \right\|_{L_{t}^{\tilde q}(0,T;L_{x}^{\tilde r}(\R^2))} \leq \epsilon^{\frac{1}{\tilde q}} \left( 1+\epsilon^{\frac14} T^{\frac34} \right) \colr \tilde C_{d}\colb,\\
&\colp \left\| \nabla^{\perp} \frac{\nabla^{\perp}}{\Delta}  \cdot \bV - \tilde \bV \right\|_{L_{t}^{\infty}(0,T;L_{x}^{2}(\R^2))} \leq \epsilon^{\frac{1}{4}} \left( 1+\epsilon^{\frac14} T^{\frac34} \right) T^{\frac34} e^{ \tilde C M  T}\tilde C_{d}\colb,\\
&\sup_{x_{0} \in \R^2} \left\| e^{-(x-x_{0})^2} \bp \zeta \\ \nabla \frac{\nabla}{\Delta} \cdot \bV \ep \right\|_{L_{t}^{2}(0,T;L^{2}_{x}(\R^2))} \leq \epsilon^{\frac12}(1+T) \colr \tilde C_{d}\colb.
\end{aligned}
\]
\et 

The proof follows from the same strategy as the proofs of Theorems \ref{rigidlid2d1} and \ref{rigidlid2d2}. One must control two new terms. Standard product estimates provide
\begin{align*}
&\left\| \mathcal{Q}[\epsilon \zeta](\bV) \right\|_{W^{4,1}_{x}} \lesssim C \left( \tfrac{1}{h_{0}} , \| \zeta \|_{H^{5}} \right) \| \bV \|_{H^{7}}^{2}\\
&\left\| -\frac{\mu}{3}((1+\epsilon \zeta)^2-1) \nabla \nabla \cdot \partial_{t} \bV - \frac{\epsilon \mu}{3} (1+\epsilon \zeta)^2 \nabla \zeta \nabla \cdot \partial_{t} \bV \right\|_{W^{4,1}_{x}} \lesssim C \left( \tfrac{1}{h_{0}} , \| \zeta \|_{H^{5}} \right) \| \epsilon \mu \partial_{t} \nabla \cdot \bV \|_{H^{6}} 
\end{align*}
and using for instance \cite[Lemmas 2.1 and 2.4]{duchene_israwi_green_naghdi} and standard product estimates we obtain
\[
\begin{aligned}
\| \epsilon \mu \partial_{t}  \nabla \cdot \bV \|_{H^{6}} &\lesssim C \left( \tfrac{1}{h_{0}} , \| \zeta \|_{H^{6}} \right) \left\| \nabla \zeta + \epsilon (\bV \cdot \nabla) \bV + \epsilon \mu \mathcal{Q}[\epsilon \zeta](\bV) \right\|_{H^{5}}\\
&\lesssim C \left( \tfrac{1}{h_{0}} , \| \zeta \|_{H^{6}} , \| \bV \|_{H^{8}} \right).
\end{aligned}
\]
We note the strategy used in the proof of Theorems \ref{rigidlid2d1} and \ref{rigidlid2d2} also provides bounds on $\| \nabla \zeta \|_{L^{2}_{t}(0,T;L^{\infty}_{x})}$. Secondly a new term also appears in the control of the rotational component. \colp We note that
\[
\nabla^{\perp} \frac{\nabla^{\perp}}{\Delta} \cdot \left( \mu \mathcal{T}[\epsilon \zeta] \partial_{t} \bV + \mu \mathcal{Q}[\epsilon \zeta](\bV) \right) = - \epsilon \mu \nabla^{\perp} \frac{\nabla^{\perp}}{\Delta} \cdot \left( h \left[ \nabla \cdot \partial_{t} \bV + (\bV \cdot \nabla) (\nabla \cdot \bV) - (\nabla \cdot \bV)^{2}  \right] \nabla \zeta \right)
\]
and using previous bounds we get
\[
\left\| \nabla^{\perp} \frac{\nabla^{\perp}}{\Delta} \cdot \left( \mu \mathcal{T}[\epsilon \zeta] \partial_{t} \bV + \mu \mathcal{Q}[\epsilon \zeta](\bV) \right)  \right\|_{L^{2}_{x}} \leq C \left( \tfrac{1}{h_{0}} , \| \zeta \|_{H^{2}} , \| \bV \|_{H^{3}} \right) \| \nabla \zeta \|_{L^{\infty}_{x}}
\]\colb
so that the strategy used in the proof of Theorems \ref{rigidlid2d1} and \ref{rigidlid2d2} to control the vorticity component can easily be adapted.

\appendix

\section{Littlewood-Paley decompostion}\label{s:LP}

In this section we introduce homogeneous and inhomogeneous Littlewood-Paley decompositions and provide basic properties. Let $\varphi_{0}$ be a smooth nonnegative even function supported in $[-1,1]$, that is equal to $1$ in $[-\frac12,\frac12]$ and that is nonincreasing on $\R^{+}$. Then we define, for any $y \in \R$ and any $j \in \mathbb{Z}$ the function $P_{j}(y) :=\varphi_{0}(2^{-j-1} y)-\varphi_{0}(2^{-j} y)$. We note that $P_{j}$ is a function supported in the annulus $\mathcal{C}( 2^{j-1}, 2^{j+1})$ for any $j \in \Z$. 

For any $y \in \R$
\[
P_{j}(y) \in [0,1] \text{   ,   } \varphi_{0}(y) + \sum_{j \in \mathbb{N}} P_{j}(y) = 1 \text{   ,   } \frac12 \leq \varphi_{0}(y)^{2} + \sum_{j \in \mathbb{N}} P^{2}_{j}(y) \leq 1.
\]
Then for any $p \in [1,\infty]$ and any Schwartz class function $f$ 
\[
\varphi_{0}(|D|) f + \sum_{j=0}^{N} P_{j}(|D|) f \overset{L^{p}}{\underset{N \to \infty}{\shortrightarrow}} f
\]
since 
\[
\| (1- \varphi_{0}(2^{-N-1} |D| )) f \|_{L^{p}} \underset{N \to \infty}{\rightarrow}  0.
\]
Such decomposition of the function $f$ is called inhomogeneous Littlewood-Paley decomposition.

For any $y \in \R^{\ast}$
\[
P_{j}(y) \in [0,1] \text{   ,   } \sum_{j \in \mathbb{Z}} P_{j}(y) = 1 \text{   ,   } \frac12 \leq \sum_{j \in \mathbb{Z}} P^{2}_{j}(y) \leq 1.
\]
Then for any $p \in (1,\infty]$ and any Schwartz class function $f$ we have
\[
\sum_{|j| \leq N} P_{j}(|D|) f \overset{L^{p}}{\underset{N \to \infty}{\shortrightarrow}} f
\]
since 
\[
\| \varphi_{0}(2^{N} |D|) f \|_{L^{p}} \underset{N \to \infty}{\rightarrow} 0 \text{   and   } \| (1- \varphi_{0}(2^{-N-1} |D| )) f \|_{L^{p}} \underset{N \to \infty}{\rightarrow}  0.
\]
Such decomposition of the function $f$ is called homogeneous Littlewood-Paley decomposition.

\section{Fourier Multipliers on Lebesgue spaces}\label{s:fouriermulti}

In this section we gather useful estimates concerning Fourier multipliers on $L^{p}$. The first lemma is about Bessel potentials.

\bl\label{control_Bessel}
Let $n=1$ or $2$. For any $\alpha \geq 0$ there exists a constant $C>0$ such that for any $a \geq 0$, any $p \in [1,\infty]$ and any $f \in L^{p}(\R^n)$
\begin{align*}
&\left\| (1+ a |D|^2)^{-\frac{\alpha}{2}} f \right\|_{L^{p}(\R^n)} \leq C \| f \|_{L^{p}(\R^n)},\\
&\left\| a^{\frac{\alpha}{2}} |D|^{\alpha} (1+ a |D|^2)^{-\frac{\alpha}{2}} f \right\|_{L^{p}(\R^n)} \leq C \| f \|_{L^{p}(\R^n)}.
\end{align*}
Furthermore for any $\alpha \geq 0$ and any $b>0$ there exists a constant $C>0$ such that for any $a \geq 0$, any $p \in [1,\infty]$ and any $f \in L^{p}(\R^n)$
\[
\left\| (1+ ba |D|^2)^{\frac{\alpha}{2}}(1+ a |D|^2)^{-\frac{\alpha}{2}} f \right\|_{L^{p}(\R^n)} \leq C \| f \|_{L^{p}(\R^n)}.
\]
Finally, for any $p \in (1,\infty)$ there exists a constant $C>0$ such that for any $f \in W^{1,p}(\R^{n})$
\[
\| |D| f \|_{L^{p}} \leq C \| \nabla f \|_{L^{p}}.
\]  
\el

\begin{proof}
By homogeneity one can assume $a=1$. As noted in \cite[V.3.1]{Stein_singular}, $\mathcal{F}^{-1}((1+ 4 \pi^2 |\xi|^2)^{-\frac{\alpha}{2}})$ is in $L^{1}(\R^{n})$ so that the first bound follows by Young's convolution inequality. The second bound is proved in \cite[V.3.2]{Stein_singular}.

Concerning the third point, we note from \cite[V.3.2]{Stein_singular} that there exists two finite measures $\nu$ and $\mu$ on $\R^{n}$ such that 
\[
(1+ ba |D|^2)^{\frac{\alpha}{2}}(1+ a |D|^2)^{-\frac{\alpha}{2}} f = \nu \ast (1+ a |D|^2)^{-\frac{\alpha}{2}} f +  \mu \ast (ba)^{\frac{\alpha}{2}} |D|^{\alpha}(1+ a |D|^2)^{-\frac{\alpha}{2}}  f
\]
so that the result follows from the first point.

Finally, since $|D| = - \sum_{j=1}^{n} \frac{\partial_{i}}{|D|} \partial_{i}$, the last point follows from the fact that the Riesz transforms are bounded on $L^{p}$ for $p \in (1,\infty)$.
\end{proof}

We then recall Bernstein's Lemma.
\bl\label{Bernstein}
Let $n \in \N^{\ast}$ and $b>a>0$. Let $\phi$ a smooth function supported in $[a,b]$ and $\chi$ a smooth function compactly supported. Then for any $s \in \R$ and any $k \in \N$ there exists a constant $C>0$ such that for any $\lambda>0$, any $p,q \in [1,\infty]$ with $q \geq p$ and any $f \in L^{p}(\R^n)$
\begin{align*}
&\left\| \nabla^{k} \chi(\lambda^{-1} |D|) f \right\|_{L^{q}(\R^{n})} \leq C \lambda^{k+d(\tfrac1p-\tfrac1q)} \left\| \chi(\lambda^{-1} |D|) f \right\|_{L^{p}(\R^{n})},\\
&\frac{1}{C} \lambda^{s} \left\| \phi(\lambda^{-1} |D|) f \right\|_{L^{p}(\R^{n})} \leq \left\| |D|^{s} \phi(\lambda^{-1} |D|) f \right\|_{L^{p}(\R^{n})} \leq C \lambda^{s} \left\| \phi(\lambda^{-1} |D|) f \right\|_{L^{p}(\R^{n})} .
\end{align*}
\el

Then we provide a high frequency result.

\bl\label{controlHFLp}
Let $\beta > 0$, $n \in \N^{\ast}$. Let $\chi$ be a smooth compactly supported function that is equal to $1$ near $0$. There exists a constant $C>0$ such that for any $p \in [1,\infty]$, any Schwartz class function $f$ and any $\lambda > 0$
\[
\left\| \frac{1- \chi( \lambda |D|)}{(\lambda |D|)^{\beta}} f \right\|_{L^{p}(\R^n)} \leq C \left\| \frac{1- \chi(\lambda |D|)}{(1+\lambda^2 |D|^2)^{\frac{\beta}{2}}} f \right\|_{L^{p}(\R^n)}.
\]
\el

\begin{proof}
By homogeneity one can assume $\lambda=1$. Using \cite[V.3.2]{Stein_singular} there exists two finite measures $\nu$ and $\mu$ on $\R^{n}$ such that 
\[
\frac{(1+ |D|^2)^{\frac{\beta}{2}}}{(1+ |D|^2)^{\frac{\beta}{2}}} \frac{1- \chi(|D|)}{|D|^{\beta}} f = \nu \ast \frac{1- \chi(|D|)}{|D|^{\beta}} \frac{1}{(1+|D|^2)^{\frac{\beta}{2}}} f +  \mu \ast\frac{1- \chi(|D|)}{(1+|D|^2)^{\frac{\beta}{2}}}  f.
\]
Then we get
\[
\left\| \mu \ast\frac{1- \chi(|D|)}{(1+|D|^2)^{\frac{\beta}{2}}}  f \right\|_{L^{p}} \lesssim \left\| \frac{1- \chi(|D|)}{(1+|D|^2)^{\frac{\beta}{2}}}  f \right\|_{L^{p}}.
\]
Furthermore, using a Littlewood-Paley decomposition as in Section \ref{s:LP} together with Bernstein's Lemma \ref{Bernstein}, there exists an integer $k_{0} \in \Z$ such that for any Schwartz class function $g$
\bas
\left\| \frac{1- \chi(|D|)}{|D|^{\beta}} g \right\|_{L^{p}} &= \left\| \sum_{j \geq k_{0}} \frac{1- \chi(|D|)}{|D|^{\beta}} P_{j}(|D|) g \right\|_{L^{p}} \leq \sum_{j \geq k_{0}} \left\|  \frac{1- \chi(|D|)}{|D|^{\beta}} P_{0}(2^{-j} |D|) g \right\|_{L^{p}}\\
&\lesssim \sum_{j \geq k_{0}} 2^{-\beta j} \left\| P_{0}(|D|) (1- \chi(|D|)) g \right\|_{L^{p}} \lesssim \| (1- \chi(|D|)) g \|_{L^{p}}
\eas
so that
\[
\left\| \nu \ast \frac{1- \chi(|D|)}{|D|^{\beta}} \frac{1}{(1+|D|^2)^{\frac{\beta}{2}}} f \right\|_{L^{p}}  \lesssim \left\| \frac{1- \chi(|D|)}{|D|^{\beta}} \frac{1}{(1+|D|^2)^{\frac{\beta}{2}}} f \right\|_{L^{p}} \lesssim \left\| \frac{1- \chi(|D|)}{(1+|D|^2)^{\frac{\beta}{2}}} f \right\|_{L^{p}}.
\]
\end{proof}

\begin{comment}
\begin{proof}
If $\chi_{1}$ is a smooth function whose support is a subset of $\{ \chi = 1 \}$ and that is equal to $1$ near $0$, we get using Mikhlin Multiplier theorem
\bas
\left\| \frac{1- \chi( \lambda |D|)}{(\lambda |D|)^{\beta}} f \right\|_{L^{p}(\R^n)} &= \lambda^{\frac{n}{p}} \left\| \frac{1- \chi_{1}(|D|)}{|D|^{\beta}} (1- \chi( |D|)) f(\lambda \cdot) \right\|_{L^{p}(\R^n)}\\
&\lesssim \lambda^{\frac{n}{p}} \left\| (1- \chi(|D|)) f(\lambda \cdot) \right\|_{L^{p}(\R^n)} = \left\| (1- \chi( \lambda |D|)) f \right\|_{L^{p}(\R^n)}.
\eas
\end{proof}
\end{comment}

In the following we provide a boundedness result in $L^{1}$ when $n=1$.
\bl\label{control|D|1d}
Let $s \in [0,1)$. There exists a constant $C_{s}>0$ such that for any $f \in W^{1,1}(\R)$
\[
\left\| |D|^{s} f \right\|_{L^{1}(\R)} \leq C_{s} \| f \|_{W^{1,1}(\R)}
\]
and there exists a constant $C>0$ such that for any $f \in W^{2,1}(\R)$
\[
\left\| |D| f \right\|_{L^{1}(\R)} \leq C \| f \|_{W^{2,1}(\R)}.
\]
\el

\begin{proof}
Using an inhomogeneous Littlewood-Paley decomposition as in Section \ref{s:LP}, we have
\[
\left\| |D|^{s} f \right\|_{L^{1}(\R)} \leq \left\| \varphi_{0}(|D|) |D|^{s} f \right\|_{L^{1}(\R)} + \sum_{j \in \N} \left\| \sgn(D) |D|^{s-1} P_{j}(|D|) \partial_{x} f \right\|_{L^{1}(\R)}.
\]
Using Lemma \ref{control_Bessel}, Young's convolution inequality and since $\xi \mapsto \varphi_{0}(|\xi|) (1+|\xi|^2)^{\frac{s}{2}}$ is a smooth compactly supported function
\[
\left\| \varphi_{0}(|D|) |D|^{s} f \right\|_{L^{1}(\R)} \lesssim \left\| \varphi_{0}(|D|) (1+|D|^2)^{\frac{s}{2}} f \right\|_{L^{1}(\R)} \lesssim \left\| f \right\|_{L^{1}(\R)}.
\]
Then we note that for any $\alpha \in \R$, the map $\xi \mapsto \sgn(\xi) |\xi|^{\alpha} P_{0}(|\xi|)$ is a smooth compactly supported function so that
\[
\left\| \mathcal{F}^{-1} ( \sgn(\xi)^{k} |\xi|^{\alpha} P_{j}(|\xi|)) \right\|_{L^{1}} = 2^{\alpha j} \left\| \mathcal{F}^{-1} ( \sgn(\xi)^{k} |\xi|^{\alpha} P_{0}(|\xi|)) \right\|_{L^{1}} \lesssim 2^{\alpha j}.
\]
Therefore it follows from Young's convolution inequality
\[
\left\| |D|^{s} f \right\|_{L^{1}(\R)} \lesssim \left\| f \right\|_{L^{1}(\R)} + \sum_{j \in \N} 2^{(s-1)j} \left\| \partial_{x} f \right\|_{L^{1}(\R)}
\]
and the first point follows. The second point follows the same way.
\end{proof}

A similar result can be obtained when $n=2$.

\bl\label{control|D|2d}
Let $s \in [0,2]$. There exists a constant $C>0$ such that for any $f \in W^{2,1}(\R^2)$
\[
\left\| |D|^{s} f \right\|_{L^{1}(\R^2)}  \leq C \| f \|_{W^{2,1}(\R^2)}.
\]
\el

\section{Dispersive estimates}\label{s:dispestim}

In this section we gather different dispersive estimates that are useful through this work. There are obtained from \cite{my_weaklydispestim}. We begin with the case $n=1$.

\bl\label{decay1d}
Let $n=1$. Let $\lambda > 0$, $\alpha \in \R$ with $\alpha \notin \{-2,-1\}$, $\beta \geq 0$ and $l \in \N$ with $l \geq 2$. Assume that $g$ is an odd $\mathcal{C}^{2}$ function. Let $y_{1}>y_{0}>0$. Let $\chi$ be a smooth even compactly supported function whose support is a subset of $[-y_{0},y_{0}]$ and that is equal to $1$ on $[-\tfrac12 y_{0},\tfrac12 y_{0}]$.

i) Let $s \in [0,\frac{\beta}{2}]$. Assume that $|g''| \geq \lambda y^{\beta}$ on $[0,y_{0}]$ and, if $s=\frac{\beta}{2}$, that $|g'-g(0)| \geq \lambda y^{\beta+1}$ on $[0,y_{0}]$. There exists $C>0$ such that for any $\mu>0$, any $t \in \R^{\ast}$, any $m \in \{ 0,1 \}$ and any Schwartz class function $f$ 
\[
\left\| e^{\iD \frac{t}{\sqrt{\mu}}  g(\sqrt{\mu} D)} \chi(\sqrt{\mu} D) (\sgn (D))^m |D|^{s} f \right\|_{L^{\infty}_{x}} \leq \frac{C}{|t|^{\frac{s+1}{2+\beta}}} \mu^{-\frac{(\beta+1)(s+1)}{2(2+\beta)}} \left\| \colr \chi(\sqrt{\mu} |D|) \colb f \right\|_{L^{1}}.
\]

ii) Assume that $g$ is $\mathcal{C}^{l}(\R)$. Assume that $\displaystyle \sum_{p=2}^{l} |g^{(p)}|\geq \lambda$ on $[\tfrac12 y_{0}, 2 y_{1}]$, that $|g''| \geq \lambda y^{\alpha}$ on $[y_{1},\infty)$ and, if $l=2$, that $\frac{1}{\lambda} y^{\alpha + 1} \geq |g' - a | \geq \lambda y^{\alpha+1}$ on $[y_{1},\infty)$ for some $a \in \R$. There exists $C>0$ such that for any $\mu>0$, any $t \in \R^{\ast}$ and any Schwartz class function $f$
\[
\left\| e^{\iD \frac{t}{\sqrt{\mu}}  g(\sqrt{\mu} D)} (1-\chi(\sqrt{\mu} D)) f \right\|_{L^{\infty}_{x}} \leq \frac{C}{|t|^{\frac{1}{l}}} \mu^{\frac{1-l}{2l}} \left\| (\sqrt{\mu} |D|)^{\frac{l-2-\alpha}{l}} (1-\chi(\sqrt{\mu} D)) f \right\|_{L^{1}}.
\]
\el

\begin{proof}
We introduce an homogeneous Littlewood-Paley decomposition as in Section \ref{s:LP}. There exists $k_{0} \in \Z$ such that using Young's convolution inequality and Bernstein's Lemma \ref{Bernstein}
\begin{align*}
&\left\| e^{\iD \frac{t}{\sqrt{\mu}}  g(\sqrt{\mu} D)} \chi(\sqrt{\mu} D) (\sgn (D))^m |D|^{s} f \right\|_{L^{\infty}_{x}}\\
&\hspace{4cm}= \left\| \sum_{k \leq k_{0}} e^{\iD \frac{t}{\sqrt{\mu}}  g(\sqrt{\mu} D)} P_{k}(\sqrt{\mu} D) (\sgn (D))^m |D|^{s} \chi(\sqrt{\mu} |D|)  f \right\|_{L^{\infty}_{x}}\\
&\hspace{4cm}\lesssim \left\| \sum_{k \leq k_{0}} \mathcal{F}^{-1} \left( e^{\iD \frac{t}{\sqrt{\mu}}  g(\sqrt{\mu} D)} P_{k}(\sqrt{\mu} D)  (\sgn (D))^m |D|^{s} \right) \right\|_{L^{\infty}_{x}} \hspace{-0.5cm} \| \colr \chi(\sqrt{\mu} |D|) \colb  f \|_{L^{1}}.
\end{align*}
The first inequality follows from \cite[Lemma 2.6]{my_weaklydispestim}.

Secondly if $\widehat{f}$ is compactly supported there exists $k_{2} , k_{1} \in \Z$ such that using Young's convolution inequality
\bas
\left\| e^{\iD \frac{t}{\sqrt{\mu}}  g(\sqrt{\mu} D)} (1-\chi(\sqrt{\mu} D)) f \right\|_{L^{\infty}_{x}} &= \left\| \sum_{k_{2} \geq k \geq k_{1}} e^{\iD \frac{t}{\sqrt{\mu}}  g(\sqrt{\mu} D)} P_{k}(\sqrt{\mu} D) (1-\chi(\sqrt{\mu} D)) f \right\|_{L^{\infty}_{x}}\\
&\hspace{-3cm}  \lesssim \left\| \sum_{k_{2} \geq k \geq k_{1}} \mathcal{F}^{-1} \left( e^{\iD \frac{t}{\sqrt{\mu}}  g(\sqrt{\mu} D)} P_{k}(\sqrt{\mu} D) |D|^{s} \right) \right\|_{L^{\infty}_{x}} \| |D|^{-s} (1-\chi(\sqrt{\mu} D)) f \|_{L^{1}}
\eas
with $s=-\frac{l-2-\alpha}{l}$. The second inequality follows from \cite[Lemma 2.6 and Lemma 2.9]{my_weaklydispestim} and by density of $\mathcal{F}^{-1}(\mathcal{C}_{c}^{\infty}(\R))$.
\end{proof}

Then we consider the case $n=2$. 

\bl\label{decay2d}
Let $n=2$. Let $\lambda > 0$, $m \in \N$, $\beta \geq 1$ and $\alpha \in \R$ with $\alpha \notin \{-2,-1\}$. Assume that $g$ is $\mathcal{C}^{3}(\R)$. Let $y_{0}>0$. Let $\chi$ be a smooth even compactly supported function whose support is a subset of $[-y_{0},y_{0}]$ and that is equal to $1$ on $[0,\tfrac12 y_{0}]$.

i) Assume that $|g''| \geq \lambda y^{\beta}$, $|g'- g'(0) | \geq \lambda y^{\beta+1}$ and $|g'| \geq \lambda$ on $[0,y_{0}]$. There exists $C>0$ such that for any $\mu>0$, any $t \in \R^{\ast}$ and any Schwartz class function $f$ 
\[
\left\| e^{\iD \frac{t}{\sqrt{\mu}}  g(\sqrt{\mu} |D|)} \chi(\sqrt{\mu} |D|) \left( \frac{\nabla}{|D|} \right)^{m} f \right\|_{L^{\infty}_{x}} \leq \frac{C}{\mu} \left( \frac{\sqrt{\mu}}{|t|} \right) ^{\frac{5+\beta}{2(2+\beta)}} \left\| \colr \chi(\sqrt{\mu} |D|) \colb  f \right\|_{L^{1}}.
\]

ii) Let $s \in \R$ such that $(s+2)(s-\alpha)<0$ or $s=\alpha$. Assume that $|g'| \geq \lambda y^{\alpha+1}$ and $|g''| \geq \lambda y^{\alpha}$ on $[\frac{y_{0}}{2},\infty)$ and, if $s=\alpha$, that $|g'| \leq \frac{1}{\lambda} y^{\alpha+1}$, $|g''| \leq \frac{1}{\lambda} y^{\alpha}$ and $|g'''| \leq \frac{1}{\lambda} y^{\alpha-1}$ on $[\frac{y_{0}}{2},\infty)$. There exists $C>0$ such that for any $\mu>0$, any $t \in \R^{\ast}$ and any Schwartz class function $f$
\[
\left\| e^{\iD \frac{t}{\sqrt{\mu}}  g(\sqrt{\mu} |D|)} (1-\chi(\sqrt{\mu} |D|)) \left( \frac{\nabla}{|D|} \right)^{m} f \right\|_{L^{\infty}_{x}} \leq \frac{C}{\mu} \left( \frac{\sqrt{\mu}}{|t|} \right)^{\frac{s+2}{2+\alpha}} \left\| (\sqrt{\mu} |D|)^{-s} (1-\chi(\sqrt{\mu} |D|)) f \right\|_{L^{1}}.
\]

iii) Assume that $\alpha<-1$. Let $s \in \R$ such that $(s+2)(s-\frac{\alpha-1}{2})<0$ or $s=\frac{\alpha-1}{2}$. Assume that $|g'| \geq \lambda$ and $|g''| \geq \lambda y^{\alpha}$ on $[\frac{y_{0}}{2},\infty)$ and, if $s=\frac{\alpha-1}{2}$, that $\frac{1}{\lambda} y^{\alpha+1} \geq |g'-a| \geq \lambda y^{\alpha+1}$ on $[\frac{y_{0}}{2},\infty)$ for some $a \in \R^{\ast}$. There exists $C>0$ such that for any $\mu>0$, any $t \in \R^{\ast}$ and any Schwartz class function $f$
\[
\left\| e^{\iD \frac{t}{\sqrt{\mu}}  g(\sqrt{\mu} |D|)} (1-\chi(\sqrt{\mu} |D|)) \left( \frac{\nabla}{|D|} \right)^{m} f \right\|_{L^{\infty}_{x}} \hspace{-0.3cm} \leq \frac{C}{\mu} \left( \frac{\sqrt{\mu}}{|t|} \right)^{\frac{2(s+2)}{3+\alpha}} \left\| (\sqrt{\mu} |D|)^{-s} (1-\chi(\sqrt{\mu} |D|)) f \right\|_{L^{1}}.
\]
\el

\begin{proof}
We begin with the case $m=0$. Point (i) follows directly from \cite[Lemma 2.12]{my_weaklydispestim}. Concerning points (ii) and (iii), by introducing a Littlewood-Paley decomposition as in Section \ref{s:LP} and proceeding as the previous lemma for high frequencies, Points (ii) and (iii) follow respectively from Lemma \cite[Lemma 2.15]{my_weaklydispestim} and \cite[Lemma 2.17]{my_weaklydispestim}. 

We now consider the case $m \geq 1$. We claim that one can easily adapt the estimates in \cite[Section 2.4]{my_weaklydispestim} to our setting. Indeed in our case one must estimate integrals under the form $\int_{\R^{+}} e^{\iD \tfrac{t}{\sqrt{\mu}} g (\sqrt{\mu} r)} \tilde J(r|x|) \chi(\sqrt{\mu} r) r dr$ or  $\int_{\R^{+}}  e^{\iD \tfrac{t}{\sqrt{\mu}} g (\sqrt{\mu} r)} \tilde J(r|x|) P(\tfrac{\sqrt{\mu} r}{2^{k}}) r^{s+1} dr$ where $\tilde J(\tau)=\int_{0}^{2\pi} e^{\iD \tau \sin(\theta)} u(\theta) d\theta$ and $u$ is a smooth periodic function ($u \equiv 1$ in \cite[Section 2.4]{my_weaklydispestim}). Similarly as $J_{0}(\tau) := \int_{0}^{2\pi} e^{\iD \tau \sin(\theta)} d\theta$, one can decompose $\tilde J$ as $\tilde J(\tau)=\tilde h_{-}(\tau) e^{\iD \tau} + \tilde h_{+}(\tau) e^{-\iD \tau}$ where, for any $p \in \N$, $|\tilde h_{\pm}^{(p)}(\tau)| \lesssim (1+|\tau|)^{-p-\frac12}$. Then, one can adapt all the results of \cite[Section 2.4]{my_weaklydispestim} replacing $J_{0}$ by $\tilde J$ so that the strategy used to prove the case $m=0$ also works.

\end{proof}

\bl\label{decay2d1/2}
Let $n=2$. Let $\lambda > 0$, $m \in \N$. Assume that $g$ is $\mathcal{C}^{2}(\R)$. Let $y_{1}>y_{0}>0$. Let $\chi$ be a smooth even compactly supported function whose support is a subset of $[-y_{0},y_{0}]$ and that is equal to $1$ on $[0,\tfrac12 y_{0}]$, $\phi$ a smooth function supported in $[\frac12,2]$ and $\tilde \chi$ a smooth compactly supported function whose support is a subset of $[y_{0},y_{1}]$

i) Assume that $|g'| \geq \lambda$ and $g''$ has a finite number of zeros on $(0,y_{0}]$. There exists $C>0$ such that for any $\mu>0$, any $t \in \R^{\ast}$, any $l \in \N$ with $l \geq 2$, any $k \in \Z$ and any Schwartz class function $f$ 
\[
\left\| e^{\iD \frac{t}{\sqrt{\mu}}  g(\sqrt{\mu} |D|)} \chi(\sqrt{\mu} |D|) \phi( 2^{-k} \sqrt{\mu} |D|) \left( \frac{\nabla}{|D|} \right)^{m} f \right\|_{L^{\infty}_{x}} \leq \frac{C}{|t|^{\frac{1}{l}}} \left( \frac{2^{k}}{\sqrt{\mu}} \right)^{2-\frac{1}{l}} \left\| \phi( 2^{-k} \sqrt{\mu} |D|) f \right\|_{L^{1}}.
\]
In particular, for any $\eta>0$ there exists $\colr C_{\eta} \colb>0$ such that for any $\mu>0$, any $t \in \R^{\ast}$, any $l \in \N$ with $l \geq 2$ and any Schwartz class function $f$ 
\[
\left\| e^{\iD \frac{t}{\sqrt{\mu}}  g(\sqrt{\mu} |D|)} \chi(\sqrt{\mu} |D|) \left( \frac{\nabla}{|D|} \right)^{m} f \right\|_{L^{\infty}_{x}} \leq \frac{\colr C_{\eta} \colb}{|t|^{\frac{1}{l}}} \left( \left\| |D|^{2-\frac{1}{l}-\eta} f \right\|_{L^{1}} + \left\| |D|^{2-\frac{1}{l}+\eta} f \right\|_{L^{1}} \right).
\]

(ii) Assume that $g$ is $\mathcal{C}^{l}(\R)$, that $\displaystyle \sum_{p=1}^{l} |g^{(p)}| \geq \lambda$ and $g''$ has a finite number of zeros on $[y_{0},y_{1} ]$. There exists $C>0$ such that for any $t \in \R^{\ast}$ and any $\mu>0$
\[
\left\| e^{\iD \frac{t}{\sqrt{\mu}}  g(\sqrt{\mu} |D|)} \tilde \chi(\sqrt{\mu} |D|) \left( \frac{\nabla}{|D|} \right)^{m} f \right\|_{L^{\infty}_{x}} \leq \frac{C}{|t|^{\frac{1}{l}}} \left\| |D|^{2-\frac{1}{l}} f \right\|_{L^{1}}.
\]

\el

\begin{proof}
The first inequality is an easy adaptation of \cite[Lemma 2.21(1)]{my_weaklydispestim}. The second bound is consequence of the first bound together with the use of a Littlewood-Paley decomposition, Bernstein's Lemma \ref{Bernstein} and the fact that
\[
\sum_{2^{k} \leq \sqrt{\mu}} \left( \frac{2^{k}}{\sqrt{\mu}} \right)^{\eta} \lesssim_{\eta} 1 \text{   ,   } \sum_{\sqrt{\mu} \leq 2^{k} \leq 2 y_{0}} \left( \frac{2^{k}}{\sqrt{\mu}} \right)^{-\eta} \leq \sum_{\sqrt{\mu} \leq 2^{k}} \left( \frac{2^{k}}{\sqrt{\mu}} \right)^{-\eta} \lesssim_{\eta} 1.
\]
The third inequality easy follows from \cite[Lemma 2.21(2)]{my_weaklydispestim}.
\end{proof}

Finally we provide Morawetz-type estimates.

\bpr\label{morawetz}
Let $n \in \N^{\ast}$ and $T>0$. Assume that $g$ is $\mathcal{C}^{1}(\R^{\ast}_{+})$. There exists $C>0$ independent of $T$ such that for any function $f \in L^{2}(\R^{n})$, any function $F$ in $L^{\infty}(0,T;L^{2}(\R^{n}))$, any $\mu>0$, any $\epsilon>0$, any $a>0$ and any $x_{0} \in \R^{n}$
\[
\int_{0}^{T} \int_{\R^{n}} \left| |g'(\sqrt{\mu} |D|)|^{\frac12} e^{\iD \frac{t}{\epsilon \sqrt{\mu}} g(\sqrt{\mu} |D|)} f \right|^2 e^{-\frac{a}{2} |x-x_{0}|^2} dx dt \leq \epsilon \frac{C}{\sqrt{a}} \| f \|_{L_{x}^{2}}^{2}.
\]
and
\[
\int_{0}^{T} \int_{\R^{n}} \left| \int_{0}^{t} |g'(\sqrt{\mu} |D|)|^{\frac12} e^{\iD \frac{(t-s)}{\epsilon \sqrt{\mu}} g(\sqrt{\mu} |D|)} F(s,\cdot) ds \right|^2 e^{-\frac{a}{2} |x-x_{0}|^2} dx dt \leq\epsilon T^2 \frac{C}{\sqrt{a}} \| F \|_{L^{2}_{s}(0,T;L_{x}^{2})}^{2}.
\]
\epr

\begin{proof}
After an appropriate change of variable in time we get from \cite[Proposition 2.28]{my_weaklydispestim}
\[
\int_{\R} \int_{\R^{n}} \left| (|g'(\sqrt{\mu} |D|)|^{\frac12} e^{\iD \frac{t}{\epsilon \sqrt{\mu}} g(\sqrt{\mu} |D|)} f)(x) \right|^2 e^{-\frac{a}{2} |x-x_{0}|^2} dx dt \leq \epsilon \frac{C}{\sqrt{a}} \| f \|_{L^{2}}^2.
\]
Then denoting by $I$ the second quantity to bound and using Jensen's inequality and the previous estimate we obtain
\[
I \leq \int_{0}^{T} T \left\| (|g'(\sqrt{\mu} |D|)|^{\frac12} e^{\iD \frac{(t-s)}{\epsilon \sqrt{\mu}} g(\sqrt{\mu} |D|)} F(s,\cdot)) e^{-\frac{a}{4} |x-x_{0}|^2} \right\|_{L^{2}_{t}(s,T;L^{2}_{x})}^2 ds \leq \epsilon T^2 \frac{C}{\sqrt{a}} \| F \|_{L^{2}_{s}(0,T;L_{x}^{2})}^{2}.
\]
\end{proof}

\bibliographystyle{alpha}
\bibliography{biblio}
\end{document}